# Predictor-Based Tracking For Neuromuscular Electrical Stimulation

**Iasson Karafyllis[*], Michael Malisoff[**], Marcio de Queiroz[***] and Miroslav Krstic[****]**

[*]Dept. of Mathematics, National Technical University of Athens, 15780, Athens, Greece, email: iasonkar@central.ntua.gr

[**]Dept. of Mathematics, Louisiana State University, Baton Rouge, LA 70803-4918, U.S.A., email: malisoff@lsu.edu

[***]Dept. of Mechanical and Industrial Eng., Louisiana State University, Baton Rouge, LA 70803, U.S.A., email: dequeiroz@me.lsu.edu

[****]Dept. of Mechanical and Aerospace Eng., University of California, San Diego, La Jolla, CA 92093-0411, U.S.A., email: krstic@ucsd.edu

**Abstract**
A new hybrid tracking controller for neuromuscular electrical stimulation is proposed. The control scheme uses sampled measurements and is designed by utilizing a numerical prediction of the state variables. The tracking error of the closed-loop system converges exponentially to zero and robustness to perturbations of the sampling schedule is exhibited. One of the novelties of our approach is the ability to satisfy a state constraint imposed by the physical system.

**Keywords:** nonlinear systems, delay systems, NMES models, numerical methods.

## 1. Introduction

Neuromuscular electrical stimulation (NMES) is a technology where skeletal muscles are artificially stimulated to help restore functionality to human limbs with motor neuron disorders [19,21]. This is done using voltage excitation of skin or implanted electrodes, which produce muscle contraction, joint torque, and limb motion. NMES is an active area of research in biomedical and rehabilitation engineering since it is key to developing neuroprosthetic devices.

To obtain a desired motion, NMES must be applied in a certain manner. NMES control is challenging due to the nonlinear, time-varying, uncertain dynamics. The problem is compounded by the presence of time delays in the muscle response, due to finite propagation of chemical ions in the muscle, synaptic transmission delays, and other causes [19]. The simplest method for generating the desired limb motion is to apply the voltage signal via open-loop control using predefined stimulation schemes specific to the functionality being restored (e.g., walking) [2]. Not surprisingly, open-loop control was found to produce unsatisfactory results [1,2,5,6]. Despite this, most NMES controllers in clinical use are open loop [6,21]. Classical feedback controllers (e.g., PID control) have also produced unsatisfactory results [16], failing to guarantee closed-loop stability [6].

In parallel to this, considerable efforts have been devoted to understanding and modeling the nonlinear physiological and mechanical dynamics of muscle stimulation, activation, and contraction [2,3,11,15]. The availability of such models has enabled researchers to explore advanced, model-based feedback control methods to improve the effectiveness of NMES. Some work along these lines includes sliding mode control [6], adaptive control [10], neural network-based controllers [4,14,21], backstepping control [16,17], and dynamic robust control [18].

While previous efforts have advanced the field of nonlinear NMES control, the issue of compensation of time delays caused by the underlying (chemical) kinetics in the NMES system has received less attention. This is an important problem due to its potential destabilizing effect on closed-loop stability [12]. Typically, the delay is modeled as an input delay to the musculoskeletal dynamics [19,20] or to the muscle activation dynamics [6]. As noted in [19], most NMES controllers have not been designed to explicitly compensate for the time delay; rather, some results have simply investigated the robustness of standard controllers to the input delay (see e.g. [6]). The first work to include time delay compensation in the design of the NMES control law was [19,20]. In these papers, PD and PID algorithms modified with a delay compensation term were designed using the predictor control approach [9]. In both, the tracking error for the knee joint angle was shown to be uniformly ultimately bounded using a Lyapunov-Krasovskii functional. Prediction uses dynamic controls to compensate arbitrarily long time delays, and therefore may sometimes work better than delay compensating controllers that have upper bounds on the allowable delays (but see [13] for non-predictive controls that compensate arbitrarily long input delays for nonlinear time-varying systems with no drift).

In this paper, we introduce a different type of predictor control for time delay compensation in the NMES system. We consider the musculoskeletal dynamics with an input delay as in [19,20], but with the constraint that the knee joint angle cannot physically exceed certain limits. Our control design is based on the hybrid, predictor feedback approach introduced in [8]. Specifically, the approach in [8] is extended to account for the nonlinear, time-varying nature of the NMES tracking control problem. The control scheme uses sampled measurements and is designed by utilizing a numerical prediction of the state variables. Our control is model based and ensures exponential tracking of the desired knee joint trajectory while satisfying the aforementioned state constraint. This is an improvement over the existing NMES results [19,20], which established the weaker ultimate boundedness condition on the tracking error under the input delay and which did not take the state constraint into account. Robustness to perturbations of the sampling schedule is also guaranteed.

The rest of this paper is organized as follows. In Section 2, we review the NMES model and make our control objectives precise. The control scheme is introduced and explained and the main result is stated (Theorem 2.1). In Section 3, we state certain general results on numerical approximation of solutions of time-varying systems, which generalize the corresponding results in [8] and use the step-size control ideas developed in [7]. The general results are used for the proof of our main tracking result, which is given in Section 4. Section 5 contains the conclusions of the present work as well as formulas for the direct application of the hybrid feedback law by the user. The Appendix contains the proofs of certain claims which are used in the proof of Theorem 2.1.

*Notation.* Throughout the paper we adopt the following notation:

* For a vector $x \in \Re^n$ we denote by $|x|$ its usual Euclidean norm, by $x'$ its transpose. The norm $|A|$ of a matrix $A \in \Re^{m \times n}$ is defined by $|A| = \max\{|Ax| : x \in \Re^n, |x| = 1\}$.

* $\Re_+$ denotes the set of non-negative real numbers. $Z_+$ denotes the set of non-negative integers. A partition $\pi = \{T_i\}_{i=0}^{\infty}$ of $\Re^+$ is an increasing sequence of times with $T_0 = 0$ and $T_i \to +\infty$. For every real $x \in \Re$, $[x]$ denotes the integer part of $x \in \Re$, i.e., $[x] = \max\{k \in Z_+ : k \le x\}$.

* We say that an increasing continuous function $\gamma : \Re^+ \to \Re^+$ is of class $K$ if $\gamma(0) = 0$. We say that an increasing continuous function $\gamma : \Re^+ \to \Re^+$ is of class $K_\infty$ if $\gamma(0) = 0$ and $\lim_{s \to +\infty} \gamma(s) = +\infty$. By $KL$ we denote the set of all continuous functions $\sigma : \Re_+ \times \Re_+ \to \Re_+$ with the properties: (i) for each $t \ge 0$ the mapping $\sigma(\cdot, t)$ is of class $K$ and (ii) for each $s \ge 0$, the mapping $\sigma(s, \cdot)$ is non-increasing with $\lim_{t \to +\infty} \sigma(s, t) = 0$.

* By $C^j(A)$ ($C^j(A; \Omega)$), where $A \subseteq \Re^n$ ($\Omega \subseteq \Re^m$), where $j \ge 0$ is a non-negative integer, we denote the class of functions (taking values in $\Omega \subseteq \Re^m$) that have continuous derivatives of order $j$ on $A \subseteq \Re^n$.

* Let $x : [a - r, b) \to \Re^n$ with $b > a \ge 0$ and $r \ge 0$. By $\breve{T}_r(t)x$ we denote the "open history" of $x$ from $t - r$ to $t$, i.e., $(\breve{T}_r(t)x)(\theta) := x(t + \theta)$ for all $\theta \in [-r, 0)$ and $t \in [a, b)$.

* Let $I \subseteq \Re^+ := [0, +\infty)$ be an interval. By $L^\infty(I; U)$ we denote the space of measurable and bounded functions $u(\cdot)$ defined on $I$ and taking values in $U \subseteq \Re^m$. We do not identify functions in $L^\infty(I; U)$ that differ on a measure zero set. For $x \in L^\infty([-r, 0); \Re^n)$ we define $\|x\|_r := \sup_{\theta \in [-r, 0)} |x(\theta)|$. Notice that $\sup_{\theta \in [-r, 0]} |x(\theta)|$ is not the essential supremum but the actual supremum and that is why the quantities $\sup_{\theta \in [-r, 0]} |x(\theta)|$ and $\sup_{\theta \in [-r, 0)} |x(\theta)|$ do not coincide in general.

* A function $h : A \to \Re$ where $0 \in A \subseteq \Re^n$ is positive definite if $h(0) = 0$ and $h(x) > 0$ for all $x \ne 0$. A function $h : \Re^n \to \Re$ is radially unbounded provided that the set $\{x \in \Re^n : h(x) \le M\}$ is bounded or empty for every $M > 0$.

## 2. The NMES Model and Main Results

We review the muscle activation and limb model from [21], which has the form

$$M_I + M_e + M_g + M_v = U$$

where $M_I$ is the inertia of the shank-foot complex about the knee joint, $M_e$ is the elasticity arising from joint stiffness, $M_g$ denotes the gravitational component, $M_v$ captures the viscous effects in the musculotendon complex from damping, and $U$ is the torque produced by the electric potential at the knee joint. The state $q$ is the angular position of the lower shank about the knee-joint. Following [21], a possible choice for the NMES model is:

$$M_I = J\ddot{q}, \quad M_g = mgl\sin(q), \quad M_v = B_1 \tanh(B_2 \dot{q}) + B_3 \dot{q} \quad \text{and} \quad M_e = k_1 q \exp(-k_2 q) + k_3 \tan(q)$$

Here $J$ and $m$ are the inertia and combined mass of the shank and foot, respectively, $l$ is the distance between the knee joint and the lumped center of the mass of the foot and shank, $k_i, B_i$ ($i=1,2,3$) are positive constants and $g$ is the gravity constant. We assume that all model parameters are known. In [21] the $k_3 \tan(q)$ term is not present, because [21] does not impose the state constraint $q(t) \in \left(-\frac{\pi}{2}, \frac{\pi}{2}\right)$ that we impose here. We will use the extra term $k_3 \tan(q)$ to ensure forward completeness of the tracking system when there is a bounded torque. The torque is the (known) total muscle force at the tendon, and has the form

$$U = \zeta(q)\eta(q,\dot{q})v$$

where $\zeta(q)$ is the known positive moment arm, $v$ is the applied voltage potential across the quadriceps muscle and $\eta(q,\dot{q})$ captures active and passive muscle characteristics and the dynamics of muscle recruitment.

We find it convenient to write the model in the form

$$\ddot{q}(t) = -\frac{dF}{dq}(q(t)) - H(\dot{q}(t)) + G(q(t),\dot{q}(t))v(t-\tau)$$
$$q(t) \in \left(-\frac{\pi}{2},\frac{\pi}{2}\right), \dot{q}(t) \in \Re, v(t) \in \Re \qquad (2.1)$$

where $F:\left(-\frac{\pi}{2},\frac{\pi}{2}\right) \to \Re_+$ is a $C^2$ non-negative function with $\lim_{q \to \pm \pi/2} F(q) = +\infty$, $H:\Re \to \Re$ is a $C^1$ function with $\inf_{x \in \Re}(xH(x)) \geq 0$, $\tau > 0$ is a constant and $G:\left(-\frac{\pi}{2},\frac{\pi}{2}\right) \times \Re \to (0,+\infty)$ is a $C^1$ positive valued and bounded function. The function $F:\left(-\frac{\pi}{2},\frac{\pi}{2}\right) \to \Re_+$ denotes the ratio of the potential and the inertia of the combined human shank-foot and machine, $H:\Re \to \Re$ denotes the ratio of the viscous torque due to damping in the musculo-tendon complex and the inertia of the combined human shank-foot and machine, $v(t)$ denotes the ratio of the delayed torque production at the knee joint and the inertia of the combined human shank-foot and machine and $\tau > 0$ is the electromechanical delay in the muscle response. Our analysis is developed for the general NMES model (2.1) under the assumptions stated above; if one wants to specify the functions $F,G,H$ with the aforementioned model characteristics then the following formulas will be useful

$$F(q) = \frac{mgl}{J}(1-\cos(q)) + \frac{k_1 \exp(-k_2 q)}{Jk_2^2}(\exp(k_2 q) - 1 - k_2 q) + \frac{k_3}{J} \ln\left(\frac{1}{\cos(q)}\right),$$
$$H(\dot{q}) = \frac{B_1}{J}\tanh(B_2 \dot{q}) + \frac{B_3}{J}\dot{q} \text{ and } G(q,\dot{q}) = \frac{1}{J}\zeta(q)\eta(q,\dot{q})$$

The control objective is the asymptotic tracking of any desired signal $q_d(t)$ that satisfies the reference system

$$\ddot{q}_d(t) = -\frac{dF}{dq}(q_d(t)) - H(\dot{q}_d(t)) + G(q_d(t),\dot{q}_d(t))v_d(t-\tau) \qquad (2.2)$$

for all $t \geq 0$, where $v_d \in C^1([-\tau,+\infty);\Re)$, and the following assumption

$$\sup_{t\geq 0}|\dot{q}_d(t)| + \sup_{t\geq 0}|v_d(t-\tau)| + \sup_{t\geq 0}|\dot{v}_d(t-\tau)| < +\infty \text{ and } \sup_{t\geq 0}|q_d(t)| < \frac{\pi}{2} \quad (2.3)$$

In this work, we extend the results provided in [8] to the time-varying case and we provide a hybrid predictor feedback controller that guarantees global asymptotic and local exponential convergence of the tracking error. Moreover, our controller does not require continuous measurement of the state variables but rather sampled measurements. The latter feature is important for practical purposes. In order to describe our obtained results, we set

$$\zeta_{1,d}(t) = \tan(q_d(t)) \text{ and } \zeta_{2,d}(t) = \frac{\dot{q}_d(t)}{\cos^2(q_d(t))}, \text{ for all } t \geq 0 \quad (2.4)$$

We also define the operator $\Omega(t_0,h,x;v) = \begin{bmatrix} \Omega_1(t_0,h,x;v) \\ \Omega_2(t_0,h,x;v) \end{bmatrix} \in \Re^2$ for all $(t_0,h,x) \in \Re_+ \times \Re_+ \times \Re^2$ and $v \in L^\infty([-\tau,+\infty);\Re)$ by means of the formulas:

$$\begin{aligned} \Omega_1(t_0,h,x;v) &:= x_1 + hx_2 \\ \Omega_2(t_0,h,x;v) &:= x_2 + \int_{t_0}^{t_0+h} g_1(\zeta_d(s)+x)ds + \int_{t_0}^{t_0+h} g_2(\zeta_d(s)+x)v(s-\tau)ds + \zeta_{2,d}(t_0) - \zeta_{2,d}(t_0+h) \end{aligned} \quad (2.5)$$

where $g_1(x) := -(1+x_1^2)\frac{dF}{dq}(\tan^{-1}(x_1)) + \frac{2x_1}{1+x_1^2}x_2^2 - (1+x_1^2)H\left(\frac{x_2}{1+x_1^2}\right)$, $g_2(x) := (1+x_1^2)G\left(\tan^{-1}(x_1), \frac{x_2}{1+x_1^2}\right)$. Let $\{T_i\}_{i=0}^\infty$ be an increasing diverging sequence with $T_0 = 0$ and $\sup_{i\geq 0}(T_{i+1}-T_i) < +\infty$. Given any initial time $t_0 \geq 0$ the sampling times will be $t_0 + T_i$ ($i=0,1,2,...$). At each sampling time $t_0 + T_i$ ($i=0,1,2,...$), we measure $(q(t_0+T_i),\dot{q}(t_0+T_i)) \in \left(-\frac{\pi}{2},\frac{\pi}{2}\right) \times \Re$. Next, we perform the following calculation:

$$\begin{aligned} z_{k+1} &= \Omega(t_0+T_i+kh_i, h_i, z_k; v) \quad k=0,...,N_i-1 \\ z_0 &= \left[\tan(q(t_0+T_i)) - \tan(q_d(t_0+T_i)), \frac{\dot{q}(t_0+T_i)}{\cos^2(q(t_0+T_i))} - \frac{\dot{q}_d(t_0+T_i)}{\cos^2(q_d(t_0+T_i))}\right]' \end{aligned} \quad (2.6)$$

where $N_i \geq 1$ is a sufficiently large integer and $h_i = \frac{\tau}{N_i}$. The preceding computations can be performed because they only require the values of $v$ on the interval $[t_0+T_i-\tau, t_0+T_i)$.

The control action $v(t)$ for $t \in [t_0+T_i, t_0+T_{i+1})$ is described by the following equations:

$$v(t) = \frac{g_2(\zeta_d(t+\tau))v_d(t) - g_1(\zeta_d(t+\tau)+\xi(t)) + g_1(\zeta_d(t+\tau)) - (1+\mu^2)\xi_1(t) - 2\mu\xi_2(t)}{g_2(\zeta_d(t+\tau)+\xi(t))} \quad (2.7)$$

where $\mu > 0$ is a constant and $\xi(t) \in \Re^2$ is given by:

$$\begin{aligned} \xi_1(t) &= e^{-\mu(t-T_i-t_0)}\left((\xi_2(T_i+t_0) + \mu\xi_1(T_i+t_0))\sin(t-T_i-t_0) + \xi_1(T_i+t_0)\cos(t-T_i-t_0)\right) \\ \xi_2(t) &= e^{-\mu(t-T_i-t_0)}\left(-(\mu\xi_2(T_i+t_0) + (1+\mu^2)\xi_1(T_i+t_0))\sin(t-T_i-t_0) + \xi_2(T_i+t_0)\cos(t-T_i-t_0)\right) \end{aligned} \quad (2.8)$$

with

$$\xi(t_0 + T_i) = z_{N_i}. \qquad (2.9)$$

The control scheme described by (2.7), (2.8) and (2.9) is a combination of:

- a numerical prediction of the error variables $x_1 = \tan(q) - \tan(q_d)$, $x_2 = \dfrac{\dot q}{\cos^2(q)} - \dfrac{\dot q_d}{\cos^2(q_d)}$ at time $t_0 + T_i + \tau$ based on the knowledge of the state variables $(q(t_0 + T_i), \dot q(t_0 + T_i)) \in \left(-\dfrac{\pi}{2}, \dfrac{\pi}{2}\right) \times \Re$ : the prediction is given by (2.9),

- an intersample prediction of the error variables $x_1 = \tan(q) - \tan(q_d)$, $x_2 = \dfrac{\dot q}{\cos^2(q)} - \dfrac{\dot q_d}{\cos^2(q_d)}$ for the time interval between two consecutive measurements: the prediction is given by (2.8), and

- the application of a nominal controller with the state variables replaced by their corresponding predicted values (predictor feedback): the control action is given by (2.7).

Our results are summarized in the following theorem.

**Theorem 2.1:** *For every $\tau > 0$, $r, \mu > 0$ and for every signal $q_d : \Re_+ \to \left(-\dfrac{\pi}{2}, \dfrac{\pi}{2}\right)$ satisfying (2.2), (2.3), there exists a locally bounded mapping $N : \Re_+ \to \{1,2,3,\ldots\}$, a constant $\omega \in \left(0, \dfrac{\mu}{2}\right)$ and a locally Lipschitz, non-decreasing function $C : \Re_+ \to \Re_+$ with $C(0) = 0$, such that for every partition $\{T_i\}_{i=0}^{\infty}$ of $\Re_+$ with $\sup_{i \geq 0}(T_{i+1} - T_i) \leq r$, for every $t_0 \geq 0$, $(q_0, \dot q_0) \in \left(-\dfrac{\pi}{2}, \dfrac{\pi}{2}\right) \times \Re$ and $v_0 \in L^{\infty}([-\tau, 0]; \Re)$, the solution $(q(t), \dot q(t), v(t)) \in \left(-\dfrac{\pi}{2}, \dfrac{\pi}{2}\right) \times \Re^2$ of the closed-loop system (2.1), (2.7), (2.8), (2.9) with*

$$N_i := N\left(\left|\left(\tan(q(t_0 + T_i)) - \tan(q_d(t_0 + T_i)), \dfrac{\dot q(t_0 + T_i)}{\cos^2(q(t_0 + T_i))} - \dfrac{\dot q_d(t_0 + T_i)}{\cos^2(q_d(t_0 + T_i))}\right)\right| + \sup_{t_0 + T_i - \tau \leq s < t_0 + T_i} |v(s) - v_d(s)|\right)$$

*and initial condition $(q(t_0), \dot q(t_0)) = (q_0, \dot q_0) \in \left(-\dfrac{\pi}{2}, \dfrac{\pi}{2}\right) \times \Re$ and $v(t_0 + s) = v_0(s)$ for $s \in [-\tau, 0)$ exists for all $t \geq t_0$ and satisfies the following inequality for all $t \geq t_0$ :*

$$\begin{aligned}&|q(t) - q_d(t)| + |\dot q(t) - \dot q_d(t)| + \sup_{t - \tau \leq s < t}|v(s) - v_d(s)| \\ &\leq \exp(-\omega(t - t_0)) C\left(\dfrac{|q_0 - q_d(t_0)| + |\dot q_0 - \dot q_d(t_0)|}{\cos^2(q_0)} + \sup_{-\tau \leq s < 0}|v_0(s) - v_d(t_0 + s)|\right)\end{aligned} \qquad (2.10)$$

Theorem 2.1 guarantees robustness to perturbations of the sampling schedule, since (2.10) holds for all sampling schedules $\{t_0 + T_i : i = 1, 2, \ldots\}$ with $\sup_{i \geq 0}(T_{i+1} - T_i) \leq r$.

# 3. Numerical Approximation of the Solutions of Time-Varying Forward Complete Systems

Consider a time-varying system of the form

$$\dot{x}(t) = f(t, x(t), u(t))$$
$$x(t) \in \Re^n, u(t) \in \Re^m, t \geq 0 \qquad (3.1)$$

under the following assumptions:

**(H1)** $f: \Re_+ \times \Re^n \times \Re^m \to \Re^n$ is a continuous vector field with $f(t,0,0) = 0$ for all $t \geq 0$, that satisfies:

$$|f(t,x,u) - f(t,y,u)| \leq L(|x|+|y|+|u|)|x-y|, \text{ for all } t \geq 0, \ x, y \in \Re^n, \ u \in \Re^m \qquad (3.2)$$

$$|f(t,x,u)| \leq (|x|+|u|)L(|x|+|u|), \text{ for all } t \geq 0, \ x \in \Re^n, \ u \in \Re^m \qquad (3.3)$$

where $L: \Re_+ \to [1,+\infty)$ is a continuous, non-decreasing function.

**(H2)** There exist a $C^2$ function $W: \Re_+ \times \Re^n \to [1,+\infty)$, a constant $c > 0$ and a function $p \in K_\infty$ such that

$$\frac{\partial W}{\partial t}(t,x) + \frac{\partial W}{\partial x}(t,x) f(t,x,u) \leq cW(t,x) + p(|u|), \text{ for all } t \geq 0, \ x \in \Re^n, \ u \in \Re^m \qquad (3.4)$$

Moreover, there exists a non-decreasing, continuous function $P: \Re_+ \to \Re_+$ such that the following inequalities hold

$$P(s) \geq 1 + \sup\left\{ \left|\frac{\partial^2 W}{\partial t^2}(t,\xi)\right| + 2sL(s)\left|\frac{\partial^2 W}{\partial t \partial x}(t,\xi)\right| + s^2 L^2(s)\left|\frac{\partial^2 W}{\partial x^2}(t,\xi)\right| : |\xi| \leq s(1+tL(s)), t \geq 0 \right\}, \text{ for all } s \geq 0 \qquad (3.5)$$

$$\left|\frac{\partial W}{\partial x}(t,x)\right| \leq \sqrt{P(|x|)}, \text{ for all } t \geq 0, \ x \in \Re^n \qquad (3.6)$$

$$|f(s,x,u) - f(t,x,u)| \leq (s-t)\sqrt{P(|x|+|u|)}, \text{ for all } t \geq 0, \ x \in \Re^n, \ u \in \Re^m \text{ and } s \geq t \qquad (3.7)$$

Finally, for every $w \geq 0$ there exists a non-decreasing, continuous function $Q_w: \Re_+ \to \Re_+$ that satisfies for all $s \geq 0$:

$$Q_w(s) \geq 1 + \sup\left\{ |x| : W(t+h,x) \leq \exp(2cw)\max_{|y|\leq s}(W(t,y)) + (2c)^{-1}\exp(2cw)p(s), \text{ for some } h \in [0,w] \text{ and } t \geq 0 \right\} \qquad (3.8)$$

Inequality (3.8) guarantees that that for each fixed $t \geq 0$ the mapping $x \to W(t,x)$ is radially unbounded (because for each $M \geq 0$ the sublevel set $\{x \in \Re^n : W(t,x) \leq M\}$ is either empty or it is contained in a ball in $\Re^n$ centered at zero with radius $Q_t(s)$, where $s = p^{-1}(2c\exp(-2ct)M)$).

The following fact is a direct consequence of assumptions (H1), (H2).

**(FACT)** *For every $\tau > 0$ there exists a function $a_\tau \in K_\infty$ such that the solution $x(t)$ of (3.1) with arbitrary initial condition $x(t_0) = x_0$, $t_0 \geq 0$ corresponding to arbitrary measurable and essentially bounded input $u:[t_0, t_0 + \tau) \to \Re^m$ satisfies*

$$|x(t)| \leq a_\tau(|x_0| + \|u\|), \text{ for all } t \in [t_0, t_0 + \tau] \qquad (3.9)$$

*where*

$$\|u\| := \operatorname*{ess\,sup}_{t \in [t_0, t_0+\tau)} |u(t)|$$

*Moreover, for every $\tau > 0$, there exists a constant $M_\tau > 0$ such that $a_\tau(s) \leq M_\tau s$ for all $s \in [0,1]$.*

More specifically, the (Fact) follows from the consideration of the solution $x(t) \in \Re^n$ of (3.1) with initial condition $x(t_0) = x_0$, $t_0 \geq 0$ corresponding to arbitrary measurable and essentially bounded input $u:[t_0, t_0 + \tau) \to \Re^m$. Inequality (3.4) implies that $W(t, x(t)) \leq \exp(c(t-t_0))W(t_0, x_0) + c^{-1}\exp(c(t-t_0))p(\|u\|)$, for all $t \geq t_0$ for which the solution exists. For $t \in [t_0, t_0 + \tau]$, we get $W(t, x(t)) \leq \exp(2c\tau)W(t_0, x_0) + (2c)^{-1}\exp(2c\tau)p(\|u\|)$ and consequently (3.8) implies that $|x(t)| \leq Q_\tau(|x_0| + \|u\|)$, for all $t \in [t_0, t_0 + \tau]$ for which the solution exists. A standard contradiction argument (based on the existence of limits from the left for maximal solutions) shows that the solution exists for all $t \in [t_0, t_0 + \tau]$. Finally, exploiting (3.3) we get the following for all $t \in [t_0, t_0 + \tau]$:

$$|x(t)| \leq |x_0| + \int_{t_0}^{t} |f(s, x(s), u(s))| ds$$

$$\leq |x_0| + \int_{t_0}^{t} L(|x(s)| + \|u\|)|x(s)|ds + \|u\|\int_{t_0}^{t} L(|x(s)| + \|u\|)ds$$

$$\leq |x_0| + L(Q_\tau(|x_0| + \|u\|) + \|u\|)\int_{t_0}^{t} |x(s)|ds + \|u\|L(Q_\tau(|x_0| + \|u\|) + \|u\|)(t - t_0)$$

Using the Gronwall-Bellman lemma we obtain for all $t \in [t_0, t_0 + \tau]$:

$$|x(t)| \leq (|x_0| + \|u\|L(Q_\tau(|x_0| + \|u\|) + \|u\|)\tau)\exp(\tau L(Q_\tau(|x_0| + \|u\|) + \|u\|))$$
$$\leq (|x_0| + \|u\|)(1 + L(Q_\tau(|x_0| + \|u\|) + |x_0| + \|u\|)\tau)\exp(\tau L(Q_\tau(|x_0| + \|u\|) + |x_0| + \|u\|))$$

which shows that (3.9) holds with $a_\tau(s) := s(1 + L(Q_\tau(s) + s)\tau)\exp(\tau L(Q_\tau(s) + s))$.

Define for all $s \geq 0$:

$$A(s) := L(Q_\tau(s) + a_\tau(s) + s) \qquad (3.10)$$

$$B(s) := L(Q_\tau(s) + a_\tau(s) + s)(a_\tau(s) + s)L(a_\tau(s) + s) \qquad (3.11)$$

Consider the following numerical scheme, which is an extension of the explicit Euler method to systems with inputs: we select a positive integer $N$ and define

$$x_{i+1} = x_i + \int_{t_0+ih}^{t_0+(i+1)h} f(s, x_i, u(s))ds, \text{ for } i = 0,\ldots, N-1 \quad (3.12)$$

for $h = \dfrac{\tau}{N}$, where $u:[t_0, t_0 + \tau) \to \Re^m$ is given.

**Theorem 3.1:** *Consider system (3.1) under Assumptions (H1)-(H2). Let $\tau > 0$ be a positive constant and let arbitrary $x_0 \in \Re^n$ $t_0 \geq 0$ and arbitrary measurable and essentially bounded input $u:[t_0, t_0 + \tau) \to \Re^m$. If $N \geq \tau \dfrac{P(Q_\tau(|x_0| + \|u\|) + \|u\|)}{c}$ then the following inequalities hold:*

$$|x(t_0 + \tau) - x_N| \leq \frac{\tau B(|x_0| + \|u\|)}{2NA(|x_0| + \|u\|)} \left( \exp(\tau A(|x_0| + \|u\|)) - 1 \right) \quad (3.13)$$

$$|x_i| \leq Q_\tau(|x_0| + \|u\|), \text{ for all } i = 0,1,\ldots, N \quad (3.14)$$

*where $x(t)$ is the solution of (3.1) with initial condition $x(t_0) = x_0$ corresponding to input $u:[t_0, t_0 + \tau) \to \Re^m$ at time $t$.*

The proof of Theorem 3.1 depends on the following technical lemmas.

**Lemma 3.2:** *Consider system (3.1) under the assumptions of Theorem 3.1. If $h \leq \dfrac{cW(t_0 + ih, x_i)}{P(|x_i| + \|u\|)}$, where $P:\Re_+ \to \Re_+$ is the function involved in (3.5), (3.6) and (3.7), then*

$$W(t_0 + (i+1)h, x_{i+1}) \leq \exp(2ch)W(t_0 + ih, x_i) + \int_{t_0+ih}^{t_0+(i+1)h} \exp(2c(t_0 + (i+1)h - s))p(|u(s)|)ds \quad (3.15)$$

**Proof of Lemma 3.2:** Define the function:

$$g(\lambda) = W(t_0 + ih + \lambda h, x_i + \lambda(x_{i+1} - x_i)) \quad (3.16)$$

for $\lambda \in [0,1]$. The following equalities hold for all $\lambda \in [0,1]$:

$$\begin{aligned}\frac{dg}{d\lambda}(\lambda) &= h\frac{\partial W}{\partial t}(t_0 + ih + \lambda h, x_i + \lambda(x_{i+1} - x_i)) + \frac{\partial W}{\partial x}(t_0 + ih + \lambda h, x_i + \lambda(x_{i+1} - x_i))(x_{i+1} - x_i) \\ \frac{d^2 g}{d\lambda^2}(\lambda) &= 2h\frac{\partial^2 W}{\partial t \partial x}(t_0 + ih + \lambda h, x_i + \lambda(x_{i+1} - x_i))(x_{i+1} - x_i) + h^2\frac{\partial^2 W}{\partial t^2}(t_0 + ih + \lambda h, x_i + \lambda(x_{i+1} - x_i)) \\ &\quad + (x_{i+1} - x_i)'\frac{\partial^2 W}{\partial x^2}(t_0 + ih + \lambda h, x_i + \lambda(x_{i+1} - x_i))(x_{i+1} - x_i)\end{aligned} \quad (3.17)$$

Moreover, notice that (3.3), (3.12) and the fact $h \leq \tau$ imply that $|x_{i+1} - x_i| \leq h(|x_i| + \|u\|)L(|x_i| + \|u\|)$ and $|x_i + \lambda(x_{i+1} - x_i)| \leq (|x_i| + \|u\|)(1 + \tau L(|x_i| + \|u\|))$. The previous inequalities in conjunction with (3.5) and (3.17) give:

$$\left| \frac{d^2 g}{d\lambda^2}(\lambda) \right| \leq h^2 P(|x_i| + \|u\|) \quad (3.18)$$

where $P:\Re_+ \to \Re_+$ is the function involved in (3.5). Furthermore, inequality (3.4) in conjunction with (3.12), (3.6), (3.7) and (3.17) gives:

$$\frac{dg}{d\lambda}(0) = h\frac{\partial W}{\partial t}(t_0+ih,x_i) + \int_{t_0+ih}^{t_0+(i+1)h} \frac{\partial W}{\partial x}(t_0+ih,x_i)f(s,x_i,u(s))ds$$

$$\leq chW(t_0+ih,x_i) + \int_{t_0+ih}^{t_0+(i+1)h} p(u(s))ds + \int_{t_0+ih}^{t_0+(i+1)h} \frac{\partial W}{\partial x}(t_0+ih,x_i)(f(s,x_i,u(s))-f(t_0+ih,x_i,u(s)))ds \quad (3.19)$$

$$\leq chW(t_0+ih,x_i) + \int_{t_0+ih}^{t_0+(i+1)h} p(u(s))ds + \frac{h^2}{2}P(|x_i|+\|u\|)$$

Combining (3.16), (3.18) and (3.19), we get:

$$W(t_0+(i+1)h,x_{i+1}) = g(1) \leq (1+ch)W(t_0+ih,x_i) + \int_{t_0+ih}^{t_0+(i+1)h} p(|u(s)|)ds + h^2 P(|x_i|+\|u\|) \quad (3.20)$$

Inequality (3.20) in conjunction with the following inequality

$$(1+ch)W(t_0+ih,x_i) + \int_{t_0+ih}^{t_0+(i+1)h} p(|u(s)|)ds + h^2 P(|x_i|+\|u\|)$$

$$\leq \exp(2ch)W(t_0+ih,x_i) + \int_{t_0+ih}^{t_0+(i+1)h} \exp(2c(t_0+(i+1)h-s))p(|u(s)|)ds$$

which holds for all $h \leq \frac{cW(t_0+ih,x_i)}{P(|x_i|+\|u\|)}$ imply that (3.15) holds. The proof is complete. ◁

**Lemma 3.3:** *Consider system (3.1) under the assumptions of Theorem 3.1. If* $h \leq \frac{c}{P(Q_\tau(|x_0|+\|u\|)+\|u\|)}$ *then*

$$W(t_0+ih,x_i) \leq \exp(2cih)W(t_0,x_0) + \int_{t_0}^{t_0+ih} \exp(2c(t_0+ih-s))p(|u(s)|)ds \text{ for all } i=0,...,N \quad (3.21)$$

*where* $Q_\tau:\Re_+ \to \Re_+$ *is the function involved in (3.8) for* $w=\tau$.

**Proof of Lemma 3.3:** The proof is by induction.

First notice that (3.21) holds for $i=0$. Suppose that it holds for some $i \in \{0,...,N-1\}$. Clearly, inequality (3.21) implies

$$W(t_0+ih,x_i) \leq \exp(2c\tau)W(t_0,x_0) + \frac{\exp(2c\tau)}{2c}p(\|u\|) \quad (3.22)$$

The previous inequality in conjunction with (3.8) implies $|x_i| \leq Q_\tau(|x_0| + \|u\|)$.

Consequently, the facts that $P: \Re_+ \to \Re_+$ is non-decreasing and $W(t_0 + ih, x_i) \geq 1$ imply $h \leq \frac{c}{P(Q_\tau(|x_0| + \|u\|) + \|u\|)} \leq \frac{cW(t_0 + ih, x_i)}{P(|x_i| + \|u\|)}$. Since $h \leq \frac{cW(t_0 + ih, x_i)}{P(|x_i| + \|u\|)}$, Lemma 3.2 shows that:

$$W(t_0 + (i+1)h, x_{i+1}) \leq \exp(2ch)W(t_0 + ih, x_i) + \int_{t_0+ih}^{t_0+(i+1)h} \exp(2c(t_0 + (i+1)h - s))p(|u(s)|)ds$$

$$\leq \exp(2ch)\left[\exp(2cih)W(t_0, x_0) + \int_{t_0}^{t_0+ih} \exp(2c(t_0 + ih - s))p(|u(s)|)ds\right] + \int_{t_0+ih}^{t_0+(i+1)h} \exp(2c(t_0 + (i+1)h - s))p(|u(s)|)ds$$

The above inequality shows that (3.21) holds with $i$ replaced by $i+1$.

The proof is complete.   ◁

**Lemma 3.4:** *Consider system (3.1) under the assumptions of Theorem 3.1. Define $e_i := x_i - x(t_0 + ih)$, $i \in \{0, ..., N\}$, where $x(t)$ is the solution of (3.1) with initial condition $x(t_0) = x_0$ corresponding to input $u: [t_0, t_0 + \tau) \to \Re^m$ and suppose that $h \leq \frac{c}{P(Q_\tau(|x_0| + \|u\|) + \|u\|)}$. Then*

$$|e_i| \leq \frac{h^2}{2} B(|x_0| + \|u\|) \frac{\exp(ihA(|x_0| + \|u\|)) - 1}{\exp(hA(|x_0| + \|u\|)) - 1}, \text{ for all } i \in \{1, ..., N\} \quad (3.23)$$

*where the functions $A$ and $B$ are defined by (3.10)-(3.11).*

**Proof of Lemma 3.4:** The following equation holds for all $i \in \{0, ..., N-1\}$, as a direct consequence of (3.12):

$$e_{i+1} = e_i + \int_{t_0+ih}^{t_0+(i+1)h} (f(s, x_i, u(s)) - f(s, x(s), u(s)))ds \quad (3.24)$$

Inequality (3.2) implies the following inequality for all $i \in \{0, ..., N-1\}$ and $s \in [t_0 + ih, t_0 + (i+1)h]$:

$$|f(s, x_i, u(s)) - f(s, x(s), u(s))| \leq L(|x_i| + |x(s)| + \|u\|)|x_i - x(s)| \quad (3.25)$$

Using the definition $e_i := x_i - x(t_0 + ih)$ and inequalities (3.3) and (3.9), we get the following inequalities for all $i \in \{0, ..., N-1\}$ and $s \in [t_0 + ih, t_0 + (i+1)h]$:

$$|x_i - x(s)| \leq |e_i| + |x(s) - x(t_0 + ih)|$$
$$\leq |e_i| + (s - t_0 - ih)\left(\max_{l \in [t_0+ih, t_0+(i+1)h]} (|x(l)|) + \|u\|\right) L\left(\max_{l \in [t_0+ih, t_0+(i+1)h]} (|x(l)|) + \|u\|\right) \quad (3.26)$$
$$\leq |e_i| + (s - t_0 - ih)(a_\tau(|x_0| + \|u\|) + \|u\|)L(a_\tau(|x_0| + \|u\|) + \|u\|)$$

Notice that all hypotheses of Lemma 3.3 hold. Therefore inequality (3.21) holds for all $i = 0, ..., N$. Recall that (3.21) implies (3.22). The previous inequality in conjunction with (3.8) implies $|x_i| \leq Q_\tau(|x_0| + \|u\|)$ for all $i = 0, \ldots, N$. Exploiting the fact that $|x_i| \leq Q_\tau(|x_0| + \|u\|)$ for all $i = 0, \ldots, N$ and (3.24), (3.25) and (3.26), we conclude that for all $i \in \{0, \ldots, N-1\}$:

$$|e_{i+1}| \leq |e_i| + hL(Q_\tau(|x_0| + \|u\|) + a_\tau(|x_0| + \|u\|) + \|u\|)e_i|$$
$$+ \frac{h^2}{2} L(Q_\tau(|x_0| + \|u\|) + a_\tau(|x_0| + \|u\|) + \|u\|)(a_\tau(|x_0| + \|u\|) + \|u\|)L(a_\tau(|x_0| + \|u\|) + \|u\|) \quad (3.27)$$
$$\leq |e_i| + hA(|x_0| + \|u\|)|e_i| + \frac{h^2}{2} B(|x_0| + \|u\|)$$

by the definitions (3.10) and (3.11) of $A$ and $B$. Then inequality (3.27) shows that the following recursive relation holds for all $i \in \{0, \ldots, N-1\}$:

$$|e_{i+1}| \leq \exp(hA(|x_0| + \|u\|))|e_i| + \frac{h^2}{2} B(|x_0| + \|u\|) \quad (3.28)$$

Using the fact $e_0 = 0$, in conjunction with relation (3.28), gives the desired inequality (3.23). The proof is complete. ◁

We are now ready to prove Theorem 3.1.

**Proof of Theorem 3.1:** All assumptions of Lemma 3.3 and Lemma 3.4 hold. Consequently, inequalities (3.21) and (3.23) hold. Inequality (3.13) follows from using the fact $\exp(hA(|x_0| + \|u\|)) - 1 \geq hA(|x_0| + \|u\|)$ and definition $h = \frac{\tau}{N}$ in conjunction with (3.23) for $i = N$. Moreover, inequality (3.21) implies (3.22). The previous inequality in conjunction with (3.8) implies (3.14). The proof is complete. ◁

Theorem 3.1 allows us to construct mappings which approximate the solution of (3.1) $\tau$ time units ahead with guaranteed accuracy level. Indeed, let $R \in C^0(\Re_+; \Re_+)$ be a positive definite function with $\liminf_{s \to 0^+} \left(\frac{R(s)}{s}\right) > 0$. Define the mapping $\Phi_{t_0} : \Re^n \times L^\infty([t_0, t_0 + \tau); \Re^m) \to \Re^n$:

$$\Phi_{t_0}(x_0, u) := x_N \quad (3.29)$$

where $x_i$, $i = 1, \ldots, N$ are defined by the numerical scheme (3.12) with $h = \frac{\tau}{N}$ and $N := N(|x_0| + \|u\|)$, where

$$N(s) := \left\lceil \tau \max\left\{\frac{a_\tau(s) + s}{2R(s)} L(a_\tau(s) + s)(\exp(\tau A(s)) - 1), \frac{P(Q_\tau(s) + s)}{c}\right\}\right\rceil + 1, \text{ for } s > 0 \quad (3.30)$$

and

$$N(0) := 1 \quad (3.31)$$

Inequality (3.13) implies that the mapping $\Phi_{t_0} : \Re^n \times L^\infty([t_0, t_0 + \tau); \Re^m) \to \Re^n$ satisfies

$$\left|\Phi_{t_0}(x_0, u) - x(t_0 + \tau)\right| \leq R(|x_0| + \|u\|) \quad (3.32)$$

Inequalities (3.13) and (3.14) in conjunction with (3.32) and (3.9) imply the following inequality:

$$\left|\Phi_{t_0}(x_0,u)\right| \leq \min\{R(|x_0|+\|u\|) + a_\tau(|x_0|+\|u\|), Q_\tau(|x_0|+\|u\|)\} \quad (3.33)$$

Notice that the mapping $N(s)$ defined by (3.30) and (3.31) is locally bounded. Indeed, there exists a constant $M_\tau > 0$ such that $a_\tau(s) \leq M_\tau s$ for all $s \geq 0$ sufficiently small. Therefore, continuity of all functions involved in (3.30) in conjunction with the fact that $\liminf\limits_{s \to 0^+}\left(\dfrac{R(s)}{s}\right) > 0$ implies that

$$\sup_{0 \leq w \leq s} N(w) < +\infty, \text{ for all } s \geq 0 \quad (3.34)$$

Therefore, we conclude:

**Corollary 3.5:** *Consider system (3.1) under the assumptions of Theorem 3.1. For every positive definite function $R \in C^0(\Re_+; \Re_+)$ with $\liminf\limits_{s \to 0^+} \dfrac{R(s)}{s} > 0$ and for every $\tau > 0$, consider the mapping $\Phi_{t_0} : \Re^n \times L^\infty([t_0, t_0+\tau); \Re^m) \to \Re^n$ defined by (3.29) for all $t_0 \geq 0$, $(x_0, u) \in \Re^n \times L^\infty([t_0, t_0+\tau); \Re^m)$, where $x_i$, $i=1,\ldots,N$ are defined by the numerical scheme (3.12) with $h = \dfrac{\tau}{N}$ and $N := N(|x_0|+\|u\|)$ is defined by (3.30), (3.31). Then inequalities (3.32), (3.33) hold for all $t_0 \geq 0$, $(x_0, u) \in \Re^n \times L^\infty([t_0, t_0+\tau); \Re^m)$, where $x(t)$ denotes the solution of (3.1) with initial condition $x(t_0) = x_0$ corresponding to input $u : [t_0, t_0+\tau) \to \Re^m$ and $\|u\| := \operatorname*{ess\,sup}\limits_{t \in [t_0, t_0+\tau)} |u(t)|$. Moreover, inequality (3.34) holds for all $s \geq 0$.*

## 4. Proof of Theorem 2.1

This section is devoted to the proof of Theorem 2.1. The proof of Theorem 2.1 is constructive and formulas will be given next for the locally bounded mapping $N : \Re_+ \to \{1,2,3,\ldots\}$ involved in the hybrid dynamic feedback law defined by (2.7), (2.8) and (2.9).

We first perform the following change of coordinates:

$$\zeta_1 = \tan(q), \quad \zeta_2 = \frac{\dot{q}}{\cos^2(q)} \quad (4.1)$$

Then $\cos^2(q) = \dfrac{1}{1+\zeta_1^2}$. Hence (4.1) and (2.4) give

$$\begin{aligned}
\dot{\zeta}_1(t) &= \zeta_2(t) \\
\dot{\zeta}_2(t) &= g_1(\zeta(t)) + g_2(\zeta(t))v(t-\tau) \\
\zeta(t) &= (\zeta_1(t), \zeta_2(t)) \in \Re^2, v(t) \in \Re
\end{aligned} \quad (4.2)$$

$$\begin{aligned}
\dot{\zeta}_{1,d}(t) &= \zeta_{2,d}(t) \\
\dot{\zeta}_{2,d}(t) &= g_1(\zeta_d(t)) + g_2(\zeta_d(t))v_d(t-\tau)
\end{aligned} \quad (4.3)$$

where $g_1(\zeta) := -(1+\zeta_1^2)\dfrac{dF}{dq}(\tan^{-1}(\zeta_1)) + \dfrac{2\zeta_1}{1+\zeta_1^2}\zeta_2^2 - (1+\zeta_1^2)H\left(\dfrac{\zeta_2}{1+\zeta_1^2}\right)$, $g_2(\zeta) := (1+\zeta_1^2)G\left(\tan^{-1}(\zeta_1), \dfrac{\zeta_2}{1+\zeta_1^2}\right)$.

Next, we define:

$$x(t) = \zeta(t) - \zeta_d(t), \quad u(t) = v(t) - v_d(t) \tag{4.4}$$

Then (4.2)-(4.3) give

$$\dot{x}(t) = f(t, x(t), u(t-\tau))$$
$$x(t) = (x_1(t), x_2(t)) \in \Re^2, u(t) \in \Re \tag{4.5}$$

where

$$f(t,x,u) := \begin{bmatrix} x_2 \\ \tilde{f}(t,x) + \tilde{g}(t,x)u \end{bmatrix}$$
$$\tilde{f}(t,x) := g_1(\zeta_d(t)+x) - g_1(\zeta_d(t)) + \left(g_2(\zeta_d(t)+x) - g_2(\zeta_d(t))\right)v_d(t-\tau) \tag{4.6}$$
$$\tilde{g}(t,x) := g_2(\zeta_d(t)+x)$$

Notice that $f(t,0,0) = 0$ for all $t \geq 0$. In order to simplify the procedure of the proof we break the proof up into three steps.

**First Step:** We show that the time-varying system (4.5) satisfies assumptions (H1)-(H2) of Section 2.

**Second Step:** Construction of $N: \Re_+ \to \{1,2,3,...\}$

**Third Step:** Rest of proof

First Step: Assumption (H1) is a direct consequence of (2.3) and definitions (4.6). Therefore, we next focus on proving assumption (H2).

Define the function:

$$W(t,x) = 1 + \frac{1}{2}\left(\frac{\zeta_{2,d}(t)+x_2}{1+(\zeta_{1,d}(t)+x_1)^2}\right)^2 + F\left(\tan^{-1}(\zeta_{1,d}(t)+x_1)\right) \tag{4.7}$$

Since $F: \left(-\frac{\pi}{2}, \frac{\pi}{2}\right) \to \Re_+$ is a $C^2$ non-negative function, it follows that $W: \Re_+ \times \Re^2 \to [1, +\infty)$ is a $C^2$ function.

The formulas (4.6) for $\tilde{f}$ and $\tilde{g}$ give the following for all $t \geq 0$, $x \in \Re^2$, $u \in \Re$:

$$\frac{\partial W}{\partial t}(t,x) + \frac{\partial W}{\partial x_1}(t,x)x_2 + \frac{\partial W}{\partial x_2}(t,x)\left(\tilde{f}(t,x) + \tilde{g}(t,x)u\right)$$
$$= -\left(\frac{\zeta_{2,d}(t)+x_2}{1+(\zeta_{1,d}(t)+x_1)^2}\right)H\left(\frac{\zeta_{2,d}(t)+x_2}{1+(\zeta_{1,d}(t)+x_1)^2}\right) \tag{4.8}$$
$$+ \left(\frac{\zeta_{2,d}(t)+x_2}{1+(\zeta_{1,d}(t)+x_1)^2}\right)G\left(\tan^{-1}(\zeta_{1,d}(t)+x_1), \frac{\zeta_{2,d}(t)+x_2}{1+(\zeta_{1,d}(t)+x_1)^2}\right)(v_d(t-\tau)+u)$$

Equation (4.8) can be seen easily by noticing that $W(t,x) = 1 + \frac{1}{2}\dot{q}^2 + F(q)$. Using the fact that $H: \Re \to \Re$ is a $C^1$ function with $\inf_{x \in \Re}(xH(x)) \geq 0$, the fact that $G: \left(-\frac{\pi}{2}, \frac{\pi}{2}\right) \times \Re \to \Re_+$ is bounded and the fact $W(t,x) \geq 1$ for all $(t,x) \in \Re_+ \times \Re^2$, we get the following from equation (4.8) for all $t \geq 0$, $x \in \Re^2$ and $u \in \Re$:

$$\frac{\partial W}{\partial t}(t,x) + \frac{\partial W}{\partial x_1}(t,x)x_2 + \frac{\partial W}{\partial x_2}(t,x)\left(\tilde{f}(t,x) + \tilde{g}(t,x)u\right) \leq$$

$$\leq \frac{1}{2}\left(\frac{\zeta_{2,d}(t) + x_2}{1 + (\zeta_{1,d}(t) + x_1)^2}\right)^2 + \frac{1}{2}\tilde{G}^2\left(v_d(t-\tau) + u\right)^2$$

$$\leq \frac{1}{2}\left(\frac{\zeta_{2,d}(t) + x_2}{1 + (\zeta_{1,d}(t) + x_1)^2}\right)^2 + \tilde{G}^2 v_d^2(t-\tau) + \tilde{G}^2 u^2$$

$$\leq \frac{1}{2}\left(\frac{\zeta_{2,d}(t) + x_2}{1 + (\zeta_{1,d}(t) + x_1)^2}\right)^2 + \tilde{G}^2 v_d^2(t-\tau)W(t,x) + \tilde{G}^2 u^2$$

where $\tilde{G} := \sup\left\{G(q,\dot{q}):(q,\dot{q}) \in \left(-\frac{\pi}{2},\frac{\pi}{2}\right) \times \Re\right\}$. The above inequality in conjunction with the definition (4.7) of $W$ implies that inequality (3.4) holds with $c := \frac{1}{2} + \tilde{G}^2 \sup_{t \geq -\tau}\left(|v_d(t)|^2\right)$ and $p(s) := \tilde{G}^2 s^2$.

The fact that there exists a continuous, non-decreasing function $P: \Re_+ \to \Re_+$ such that inequalities (3.5), (3.6) and (3.7) hold is shown in the Appendix.

We next turn to (3.8). Since $F: \left(-\frac{\pi}{2}, \frac{\pi}{2}\right) \to \Re_+$ is a smooth, non-negative function with $\lim_{q \to \pm \pi/2} F(q) = +\infty$, it follows that the function $\Re^2 \ni x \to \tilde{W}(x) = 1 + \frac{1}{2}\left(\frac{x_2}{1+x_1^2}\right)^2 + F(\tan^{-1}(x_1))$ is a smooth, positive valued, radially unbounded function. Consequently, there exists a pair of $K_\infty$ functions $\theta_i: \Re_+ \to \Re_+$ ($i=1,2$) and a constant $R_2 \geq 0$ such that

$$\theta_1(|x|) \leq \tilde{W}(x) \leq R_2 + \theta_2(|x|), \text{ for all } x \in \Re^2 \tag{4.9}$$

Notice that the identity $W(t,x) = \tilde{W}(\zeta_d(t) + x)$ holds for all $(t,x) \in \Re_+ \times \Re^2$.

Therefore, it follows from (4.9) that the inequalities

$$|x| \leq \theta_1^{-1}(W(t,x)) + \sup_{l \geq 0}(|\zeta_d(l)|) \text{ and } W(t,x) \leq R_2 + \theta_2\left(\sup_{l \geq 0}(|\zeta_d(l)|) + |x|\right) \tag{4.10}$$

hold for all $(t,x) \in \Re_+ \times \Re^2$. Inequalities (4.10) imply that (3.8) holds with

$$Q_w(s) := 1 + \theta_1^{-1}\left(\exp(2cw)\left(R_2 + \theta_2\left(s + \sup_{l \geq 0}(|\zeta_d(l)|)\right)\right) + (2c)^{-1}\exp(2cw)p(s)\right) + \sup_{l \geq 0}(|\zeta_d(l)|) \tag{4.11}$$

More specifically, for each $t \geq 0$ and $s \geq 0$ (4.10) gives $\max_{|y| \leq s} W(t, y) \leq R_2 + \theta_2 \left( \sup_{l \geq 0} (|\zeta_d(l)|) + s \right)$ and consequently, the condition $W(t+h, x) \leq \exp(2cw) \max_{|y| \leq s}(W(t, y)) + (2c)^{-1} \exp(2cw) p(s)$ implies that $W(t+h, x) \leq \exp(2cw) \left( R_2 + \theta_2 \left( s + \sup_{l \geq 0} (|\zeta_d(l)|) \right) \right) + (2c)^{-1} \exp(2cw) p(s)$ for all $h \in [0, w]$. The assertion follows from the first inequality in (4.10) and the previous inequality.

Second Step: Define for all $(t, x) \in \mathfrak{R}_+ \times \mathfrak{R}^2$:

$$V(x) = \frac{2}{\mu^2 + 2 - \mu\sqrt{\mu^2 + 4}} \left( x_1^2 + (x_2 + \mu x_1)^2 \right) \tag{4.12}$$

$$k(t, x) := -\frac{(1 + \mu^2) x_1 + 2\mu x_2 + \tilde{f}(t, x)}{\tilde{g}(t, x)} \tag{4.13}$$

Definitions (4.12) and (4.13) allow us to conclude that the following relations hold:

$$\frac{\partial V}{\partial x_1}(x) x_2 + \frac{\partial V}{\partial x_2}(x) \left( \tilde{f}(t, x) + \tilde{g}(t, x) k(t, x) \right) = -2\mu V(x), \quad \forall (t, x) \in \mathfrak{R}_+ \times \mathfrak{R}^2 \tag{4.14}$$

$$|x|^2 \leq V(x) \leq K |x|^2, \quad \forall x \in \mathfrak{R}^2 \tag{4.15}$$

$$|\nabla V(x)| \leq 2K |x|, \quad \forall x \in \mathfrak{R}^2 \tag{4.16}$$

where $K := \dfrac{\mu^2 + 2 + \mu\sqrt{\mu^2 + 4}}{\mu^2 + 2 - \mu\sqrt{\mu^2 + 4}}$.

Definitions (4.6) and (4.13) allow us to assume the existence of a function $\tilde{a} \in K_\infty$, a continuous, non-decreasing function $M : \mathfrak{R}_+ \to [1, +\infty)$ and positive constants $\tilde{k}$ and $\varepsilon$ that satisfy:

$$|k(t, x)| \leq \tilde{a}(|x|), \text{ for all } (t, x) \in \mathfrak{R}_+ \times \mathfrak{R}^2 \tag{4.17}$$

$$\tilde{a}(s) := \tilde{k} s, \text{ for all } s \in [0, \varepsilon] \tag{4.18}$$

$$\tilde{g}(t, x) |k(t, x) - k(t, \xi)| \leq M(|x| + |\xi|) |\xi - x|, \text{ for all } z, x \in \mathfrak{R}^n \tag{4.19}$$

The existence of $\tilde{a} \in K_\infty$ follows from the boundedness of $\zeta_d(t)$ (see (2.3), (2.4)). The existence of $M : \mathfrak{R}_+ \to [1, +\infty)$ is shown in the Appendix. Moreover, define for all $s \geq 0$:

$$D_r(s) := 2K(a_r(s) + s) M(a_r(s) + s) \exp(rL(a_r(s) + s)), \quad \beta(s) := \tilde{a}(s\sqrt{K}) + s\sqrt{K} \tag{4.20}$$

where $a_r \in K_\infty$ is the function involved in (3.9) for system (4.5) with $\tau$ replaced by $r > 0$ and $L : \mathfrak{R}_+ \to [1, +\infty)$ is the function involved in assumption (H1) for the vector field $f(t, x, u) := \begin{bmatrix} x_2 \\ \tilde{f}(t, x) + \tilde{g}(t, x) u \end{bmatrix}$ (from the right hand side of system (4.5)).

Next select a constant $\delta > 0$, such that:

$$2\sqrt{K\delta} \leq \varepsilon \tag{4.21}$$

Having selected $\delta > 0$, we are in a position to select a constant $\gamma > 0$, so that:

$$\gamma \leq \min\left\{\sqrt{\delta}, \frac{\mu\delta}{2}\right\} \tag{4.22}$$

Define:

$$\phi := 2KM\left(2\sqrt{K\delta} + \sqrt{\delta}\right)\exp(r\tilde{L}) \tag{4.23}$$

$$\tilde{L} := L\left((1+\tilde{k})2\sqrt{K\delta} + \sqrt{\delta}\right) \tag{4.24}$$

and moreover, select a constant $\tilde{R} > 0$, so that:

$$\tilde{R}\tilde{k}\sqrt{K} < 1 \text{ and } \frac{\phi\tilde{R}}{\mu\sqrt{2}}\left(1 + \frac{\tilde{k}\sqrt{K}(\tilde{R}+1)}{1-\tilde{R}\tilde{k}\sqrt{K}}\right) < 1 \tag{4.25}$$

Finally, define:

$$R(s) := \min\left\{\frac{\gamma}{\max\{1, D_r(a_\tau(s) + \beta(Q_\tau(s)))\}}, \tilde{R}s, \frac{1}{2\sqrt{K}}\tilde{a}^{-1}\left(\frac{s}{2}\right)\right\} \tag{4.26}$$

Equation (4.18), definition (4.26) of $R$ and the fact that $Q_\tau(s) \geq 1$ for all $s \geq 0$, imply that $\liminf_{s\to 0^+}\left(\frac{R(s)}{s}\right) = \min\left\{\tilde{R}, \frac{1}{4\tilde{k}\sqrt{K}}\right\} > 0$. Therefore, Corollary 3.5 guarantees that the mapping $N: \Re_+ \to \{1,2,3,...\}$ defined by (3.30)-(3.31) for system $\dot{x} = f(t,x,u)$ (i.e., the delay-free version of (4.5)) is locally bounded and the mapping $\Phi_{t_0}: \Re^2 \times L^\infty([t_0, t_0+\tau); \Re) \to \Re^2$ defined by (3.29) satisfies inequalities (3.32)-(3.33) for all $(t_0, x_0) \in \Re_+ \times \Re^2$ and $u \in L^\infty([t_0, t_0+\tau); \Re)$, where $x(t)$ denotes the solution of $\dot{x} = f(t,x,u)$, initial condition $x(t_0) = x_0$ corresponding to input $u: [t_0, t_0+\tau) \to \Re$ and $\|u\| := \underset{t \in [t_0, t_0+\tau)}{\text{ess sup}} |u(t)|$.

In the new coordinate system defined by (4.1), (2.4), (4.4), the closed-loop system (2.1), (2.7), (2.8), (2.9) with $N_i := N\left(\left|\left(\tan(q(t_0+T_i)) - \tan(q_d(t_0+T_i)), \frac{\dot{q}(t_0+T_i)}{\cos^2(q(t_0+T_i))} - \frac{\dot{q}_d(t_0+T_i)}{\cos^2(q_d(t_0+T_i))}\right)\right| + \underset{t_0+T_i-\tau \leq s < t_0+T_i}{\sup}|v(s) - v_d(s)|\right)$

is described by the equations:

$$\dot{x}(t) = f(t, x(t), u(t-\tau))$$
$$x(t) \in \Re^2, u(t) \in \Re \tag{4.27}$$

with

$$\dot{\xi}(t) = f(t+\tau, \xi(t), k(t+\tau, \xi(t))), \xi(t) \in \Re^2$$
$$u(t) = k(t+\tau, \xi(t))$$
, for $t \in [t_0+T_i, t_0+T_{i+1})$ \tag{4.28}

and

$$\xi(t_0+T_i) = z_{N_i} \tag{4.29}$$

where $N_i := N\left(|x(t_0+T_i)| + \underset{t_0+T_i-\tau \leq s < t_0+T_i}{\sup}|u(s)|\right)$, $h_i = \frac{\tau}{N_i}$ and

$$z_{j+1} = z_j + \int_{t_0+T_i+jh_i}^{t_0+T_i+(j+1)h_i} f(s, z_j, u(s-\tau))ds, \text{ for } j = 0,\ldots, N_i - 1 \text{ and } z_0 = x(t_0+T_i) \quad (4.30)$$

To verify that (2.8) and (4.28) agree, notice that (4.13) and the fact that $f(t,x,u) := \begin{bmatrix} x_2 \\ \tilde{f}(t,x) + \tilde{g}(t,x)u \end{bmatrix}$ imply that (4.28) can be written as

$$\begin{aligned} \dot{\xi}_1(t) &= \xi_2(t) \\ \dot{\xi}_2(t) &= -(1+\mu^2)\xi_1(t) - 2\mu\xi_2(t) \end{aligned}, \text{ for } t \in [t_0+T_i, t_0+T_{i+1})$$

The solution of the above system of differential equations is given by (2.8).

<u>Third Step:</u> It should be emphasized that for every $t_0 \geq 0$, $(q_0, \dot{q}_0) \in \left(-\frac{\pi}{2}, \frac{\pi}{2}\right) \times \Re$ and $v_0 \in L^\infty([-\tau, 0); \Re)$, the solution $(q(t), \dot{q}(t), v(t)) \in \left(-\frac{\pi}{2}, \frac{\pi}{2}\right) \times \Re^2$ of the closed-loop system (2.1), (2.7), (2.8), (2.9) with $N_i := N\left(\left|\left(\tan(q(t_0+T_i)) - \tan(q_d(t_0+T_i)), \frac{\dot{q}(t_0+T_i)}{\cos^2(q(t_0+T_i))} - \frac{\dot{q}_d(t_0+T_i)}{\cos^2(q_d(t_0+T_i))}\right)\right| + \sup_{t_0+T_i-\tau \leq s < t_0+T_i}|v(s) - v_d(s)|\right)$ and initial condition $(q(t_0), \dot{q}(t_0)) = (q_0, \dot{q}_0) \in \left(-\frac{\pi}{2}, \frac{\pi}{2}\right) \times \Re$ and $v(t_0+s) = v_0(s)$ for $s \in [-\tau, 0)$ is related to the solution $(x(t), u(t)) \in \Re^3$ of the closed-loop system (4.27), (4.28), (4.29), (4.30) with initial condition $x(t_0) = \left(\tan(q(t_0)) - \tan(q_d(t_0)), \frac{\dot{q}(t_0)}{\cos^2(q(t_0))} - \frac{\dot{q}_d(t_0)}{\cos^2(q_d(t_0))}\right)$, $u(t_0+s) = v_0(s) - v_d(t_0+s)$ by means of the equations:

$$x(t) = \left(\tan(q(t)) - \tan(q_d(t)), \frac{\dot{q}(t)}{\cos^2(q(t))} - \frac{\dot{q}_d(t)}{\cos^2(q_d(t))}\right), \ u(t-\tau) = v(t-\tau) - v_d(t-\tau) \quad (4.31)$$

$$\begin{aligned} q(t) &= \tan^{-1}(x_1(t) + \tan(q_d(t))) \\ \dot{q}(t) &= \frac{1}{1+(x_1(t)+\tan(q_d(t)))^2}\left(x_2(t) + \frac{\dot{q}_d(t)}{\cos^2(q_d(t))}\right) \end{aligned} \quad (4.32)$$

which hold for all $t \geq t_0$ for which the solutions exist. The global relations

$$|q(t) - q_d(t)| \leq |x_1(t)|, \ |\dot{q}(t) - \dot{q}_d(t)| \leq |x_2(t)| + M_1|x_1(t)|,$$

$$|x_1(t_0)| \leq M_2 \frac{|q(t_0) - q_d(t_0)|}{|\cos(q(t_0))|}, \ |x_2(t_0)| \leq \frac{|\dot{q}(t_0) - \dot{q}_d(t_0)| + 2M_2M_3^2|q(t_0) - q_d(t_0)|}{\cos^2(q(t_0))},$$

where $M_1 := 2\sup_{t \geq 0}(|\zeta_{2,d}(t)|)$, $M_2 := \sup_{t \geq 0}\left(\frac{1}{|\cos(q_d(t))|}\right)$ and $M_3 := \sup_{t \geq 0}(|\dot{q}_d(t)|)$, (which are direct consequences of (4.31)-(4.32) and the Mean Value Theorem), allow us to conclude that in order to prove Theorem 2.1, it suffices to show that there exists a locally Lipschitz, non-decreasing function $\hat{C}: \Re_+ \to \Re_+$ with $\hat{C}(0) = 0$, such that for every partition $\{T_i\}_{i=0}^\infty$ of $\Re_+$ with $\sup_{i \geq 0}(T_{i+1} - T_i) \leq r$, and every $t_0 \geq 0$, $x_0 \in \Re^2$ and $u_0 \in L^\infty([-\tau, 0); \Re)$, the solution $(x(t), u(t)) \in \Re^3$ of the closed-loop system given by (4.27), (4.28), (4.29) and (4.30) with initial condition $x(t_0) = x_0$ and $u(t_0 + s) = u_0(s)$ for $s \in [-\tau, 0)$ exists for all $t \geq t_0$ and satisfies the following inequality for all $t \geq t_0$:

$$|x(t)| + \sup_{t-\tau \leq s < t} |u(s)| \leq \exp(-\omega(t-t_0))\hat{C}\left(|x_0| + \sup_{-\tau \leq s < 0} |u_0(s)|\right) \tag{4.33}$$

Therefore, the rest of proof is devoted to the proof of estimate (4.33) for the closed-loop system (4.27), (4.28), (4.29), (4.30).

Having completed the design of the feedback law by constructing the function $N : \Re_+ \to \{1,2,3,...\}$ in (3.30), we are now ready to prove some basic results concerning the closed-loop system (4.27), (4.28), (4.29), (4.30).

The following claim shows that practical stabilization is achieved. Its proof is provided in the Appendix.

**Claim 1:** *There exists $\sigma \in KL$ such that for every partition $\{T_i\}_{i=0}^{\infty}$ of $\Re_+$ with $\sup_{i \geq 0}(T_{i+1} - T_i) \leq r$, for every $(t_0, x_0) \in \Re_+ \times \Re^2$ and $u_0 \in L^{\infty}([-\tau, 0); \Re)$, the solution of (4.27), (4.28), (4.29) and (4.30) with initial condition $x(t_0) = x_0$, $\breve{T}_\tau(t_0)u = u_0$ satisfies the following inequality for all $t \geq t_0$:*

$$V(x(t)) \leq \max\left\{\sigma\left(|x_0| + \|u_0\|_\tau, t - t_0\right), \mu^{-1}\gamma\right\} \tag{4.34}$$

*where $\gamma > 0$ is the constant involved in (4.22) and (4.26).*

The following claim shows that local exponential stabilization is achieved. Its proof is provided in the Appendix.

**Claim 2:** *There exist positive constants $S_1, S_2$ and $\omega$ such that for each partition $\{T_i\}_{i=0}^{\infty}$ of $\Re_+$ with $\sup_{i \geq 0}(T_{i+1} - T_i) \leq r$ and each $(t_0, x_0) \in \Re_+ \times \Re^2$ and $u_0 \in L^{\infty}([-\tau, 0); \Re)$, the solution of (4.27), (4.28), (4.29) and (4.30) with initial conditions $x(t_0) = x_0$ and $\breve{T}_\tau(t_0)u = u_0$ satisfies the following inequalities:*

$$|u(t)|\exp(\omega(t - t_0 - T_j)) \leq S_1\left(\sup_{t_0 + T_j \leq w \leq t_0 + T_j + \tau}(|x(w)|) + \|\breve{T}_\tau(t_0 + T_j)u\|_\tau\right), \text{ for all } t \geq t_0 + T_j \tag{4.35}$$

$$|x(t)|\exp(\omega(t - t_0 - T_j - \tau)) \leq S_2\left(\sup_{t_0 + T_j \leq w \leq t_0 + T_j + \tau}(|x(w)|) + \|\breve{T}_\tau(t_0 + T_j)u\|_\tau\right), \text{ for all } t \geq t_0 + T_j + \tau \tag{4.36}$$

*where $j$ is the smallest integer for which it holds $V(x(t_0 + T_j + \tau)) \leq \delta$ and $\delta > 0$ is the constant involved in (4.21) and (4.22).*

The following claim guarantees that $u$ is bounded. Its proof is provided in the Appendix.

**Claim 3:** *There exists a non-decreasing function $S : \Re_+ \to \Re_+$ such that for each partition $\{T_i\}_{i=0}^{\infty}$ of $\Re_+$ with $\sup_{i \geq 0}(T_{i+1} - T_i) \leq r$ and each $(t_0, x_0) \in \Re_+ \times \Re^2$ and $u_0 \in L^{\infty}([-\tau, 0); \Re)$, the solution of (4.27), (4.28), (4.29) and (4.30) with initial condition $x(t_0) = x_0$ and $\breve{T}_\tau(t_0)u = u_0$ satisfies the following inequality for all $t \geq t_0$:*

$$|x(t)| + \|\breve{T}_\tau(t)u\|_\tau \leq S(|x_0| + \|u_0\|_\tau) \tag{4.37}$$

We are now ready to prove estimate (4.33). Let an arbitrary partition $\{T_i\}_{i=0}^{\infty}$ of $\Re_+$ with $\sup_{i\geq 0}(T_{i+1}-T_i)\leq r$, $(t_0,x_0)\in\Re_+\times\Re^2$ and $u_0\in L^{\infty}([-\tau,0);\Re)$ be given and consider the solution of (4.27), (4.28), (4.29) and (4.30) with initial condition $x(t_0)=x_0$ and $\breve{T}_\tau(t_0)u=u_0$.

Inequalities (4.15) and (3.9) imply that the smallest integer $j$ for which $V(x(t_0+T_j+\tau))\leq\delta$ holds is $j=0$ for the case $K(a_\tau(|x_0|+\|u_0\|_\tau))^2\leq\delta$. Moreover, the fact that there exists a constant $M_\tau>0$ such that $a_\tau(s)\leq M_\tau s$ for all $s\in[0,1]$, in conjunction with inequalities (3.9), (4.35) and (4.36), allow us to conclude that that there exists a constant $\tilde{\Omega}>0$ such that

$$|x(t)|+\|\breve{T}_\tau(t)u\|_\tau \leq \tilde{\Omega}\exp(-\omega(t-t_0))(|x_0|+\|u_0\|_\tau), \text{ for all } t\geq t_0 \tag{4.38}$$

provided that $|x_0|+\|u_0\|_\tau \leq \min\left\{1,\frac{1}{M_\tau}\sqrt{\frac{\delta}{K}}\right\}$.

Proposition 7 in [22] (which provides upper bounds for $KL$ functions $\sigma$ of the form $\sigma(s,t)\leq\beta_1(\exp(-t)\beta_2(s))$ for functions $\beta_1,\beta_2\in K_\infty$) in conjunction with (4.34), (4.22) and the fact that $\sup_{i\geq 0}(T_{i+1}-T_i)\leq r$, allow us to guarantee the existence of a non-decreasing function $\tilde{T}:\Re_+\to\Re_+$ such that the smallest sampling time $t_0+T_j$ for which $V(x(t_0+T_j+\tau))\leq\delta$ holds satisfies $T_j\leq\tilde{T}(|x_0|+\|u_0\|_\tau)$ for all $(t_0,x_0)\in\Re_+\times\Re^2$ and $u_0\in L^{\infty}([-\tau,0);\Re)$. Combining (4.35), (4.36), (4.37) with the previous inequality, allows us to conclude the existence of a non-decreasing function $\tilde{G}:\Re_+\to\Re_+$ such that the following inequality holds for all $t\geq t_0$:

$$|x(t)|+\|\breve{T}_\tau(t)u\|_\tau \leq \exp(-\omega(t-t_0))\tilde{G}(|x_0|+\|u_0\|_\tau) \tag{4.39}$$

Consequently, using (4.38) and (4.39) we conclude that (4.33) holds with $\hat{C}(s):=\frac{1}{s}\int_s^{2s}\tilde{C}(w)dw$ for all $s>0$ and $\hat{C}(0):=0$, where $\tilde{C}(s):=\max\left\{1,\frac{\tilde{G}(l)}{\tilde{\Omega}l}\right\}\tilde{\Omega}s$, for all $s\in[0,l]$, $\tilde{C}(s):=\max\{\tilde{\Omega}s,\tilde{G}(s)\}$, for all $s>l$, $l:=\min\left\{1,\frac{1}{M_\tau}\sqrt{\frac{\delta}{K}}\right\}$. The proof of Theorem 2.1 is complete. ◁

## 5. Concluding Remarks

A hybrid tracking controller for neuromuscular electrical stimulation was proposed. The main advantages of the proposed control scheme are:

- the control scheme uses sampled measurements and does not require continuous measurements of the state variables,
- the tracking error of the closed-loop system converges exponentially to zero for all initial conditions,
- the controller is designed in such a way a specific state constraint imposed by the physical system is satisfied, and
- robustness to perturbations of the sampling schedule is guaranteed.

The control scheme (2.7), (2.8), (2.9) can be programmed easily. However, it requires knowledge of the signal to be tracked $q_d(t)$, the specific functions $F, G$ and $H$, the delay $\tau$ appearing in the NMES model (2.1) and the upper diameter of the sampling schedule $r$; it is a model-based nonlinear hybrid predictor feedback. Knowledge of the aforementioned functions and constants can lead the user to an easy implementation of the proposed control scheme by utilizing the formulas in Tables 1 and 2: both tables contain all the formulas and constants which are involved in the control scheme and are selected in such a way that all inequalities and equalities used in previous sections are satisfied automatically. The parameters $\mu > 0$ and $\varepsilon > 0$ are the controller parameters (to be selected by the user).

However, the fact that the proposed control scheme is model-based is possibly a disadvantage; the robustness with respect to modelling errors of the NMES model has to be studied.

| | |
|---|---|
| $c$ | $= \dfrac{1}{2} + \tilde{G}^2 \Lambda_2^2$ |
| $\tilde{G}$ | $= \sup\left\{ G(q,\dot{q}) : (q,\dot{q}) \in \left(-\dfrac{\pi}{2}, \dfrac{\pi}{2}\right) \times \Re \right\}$ |
| $K$ | $= \dfrac{\mu^2 + 2 + \mu\sqrt{\mu^2 + 4}}{\mu^2 + 2 - \mu\sqrt{\mu^2 + 4}}$ |
| $\gamma$ | $= \min\left\{ \dfrac{\varepsilon}{2\sqrt{K}}, \dfrac{\mu\varepsilon^2}{8K} \right\}$ |
| $\phi$ | $= 2KM\left(\varepsilon + \dfrac{\varepsilon}{2\sqrt{K}}\right)\exp(r\tilde{L})$ |
| $\tilde{L}$ | $= L\left((1+\tilde{k})\varepsilon + \dfrac{\varepsilon}{2\sqrt{K}}\right)$ |
| $\tilde{R}$ | $= \min\left\{ \dfrac{\mu\sqrt{2}}{2\phi\left(1+4\tilde{k}\sqrt{K}\right)}, \dfrac{1}{2\tilde{k}\sqrt{K}}, \dfrac{1}{2} \right\}$ |
| $\Lambda_1$ | $= \sup_{t \geq 0}\left|\zeta_d(t)\right|$ |
| $\Lambda_2$ | $= \sup_{t \geq -\tau}\left|v_d(t)\right|$ |
| $\Lambda_3$ | $= \sup_{t \geq 0}\left|\dot{\zeta}_d(t)\right|$ |
| $\Lambda_4$ | $= \sup_{t \geq -\tau}\left|\dot{v}_d(t)\right|$ |
| $\Lambda_5$ | $= \sup_{t \geq 0}\left|\ddot{\zeta}_d(t)\right|$ |
| $\tilde{k}$ | $= \dfrac{(1+\mu)^2 + \psi_1(\varepsilon) + \Lambda_2\psi_2(\varepsilon)}{\min\left\{g_2(\zeta): |\zeta| \leq \Lambda_1 + \varepsilon\right\}}$ |

**Table 1:** Table of all constants involved the hybrid feedback law (2.6), (2.7), (2.8), (2.9). The functions $\psi_i(s)$ $(i = 1,2)$ are arbitrary continuous, non-decreasing functions that satisfy $\max\left\{|\nabla g_i(\zeta)| : |\zeta| \leq \Lambda_1 + s\right\} \leq \psi_i(s)$, $i = 1,2$ for all $s \geq 0$.

| | |
|---|---|
| $N(s)$ | $= \left[\tau \max\left\{\frac{a_\tau(s)+s}{2R(s)} L(a_\tau(s)+s)(\exp(\tau A(s))-1), c^{-1}P(Q_\tau(s)+s)\right\}\right]+1$, if $s>0$<br>$=1$, if $s=0$ |
| $L(s)$ | $=1+s+\psi_1(s)+(1+\Lambda_2)\psi_2(s)+\widetilde{G}(1+2\Lambda_1^2+2s^2)$ |
| $(g_1(\zeta), g_2(\zeta))$ | $=\left(-(1+\zeta_1^2)\frac{dF}{dq}(\tan^{-1}(\zeta_1))+\frac{2\zeta_1}{1+\zeta_1^2}\zeta_2^2-(1+\zeta_1^2)H\left(\frac{\zeta_2}{1+\zeta_1^2}\right), (1+\zeta_1^2)G\left(\tan^{-1}(\zeta_1), \frac{\zeta_2}{1+\zeta_1^2}\right)\right)$ |
| $A(s)$ | $=L(Q_\tau(s)+a_\tau(s)+s)$ |
| $\widetilde{W}(x)$ | $=1+\frac{1}{2}\left(\frac{x_2}{1+x_1^2}\right)^2+F(\tan^{-1}(x_1))$ |
| $a_\tau(s)$ | $=s(1+L(Q_\tau(s)+s)\tau)\exp(\tau L(Q_\tau(s)+s))$ |
| $Q_\tau(s)$ | $=1+\theta_1^{-1}\left(\exp(2c\tau)(R_2+\theta_2(s+\Lambda_1))+(2c)^{-1}\exp(2c\tau)\widetilde{G}^2 s^2\right)+\Lambda_1$ |
| $R(s)$ | $=\min\left\{\frac{\gamma}{\max\{1, D_r(a_\tau(s)+\beta(Q_\tau(s)))\}}, \widetilde{R}s, \frac{1}{2\sqrt{K}}\widetilde{a}^{-1}\left(\frac{s}{2}\right)\right\}$ |
| $D_r(s)$ | $=2K(a_r(s)+s)M(a_r(s)+s)\exp(rL(a_r(s)+s))$ |
| $M(s)$ | $=\left((1+\mu)^2+L(s)\right)\left(1+\frac{s\psi_2(s)}{\min\{g_2(\zeta):|\zeta|\leq\Lambda_1+s\}}\right)$ |
| $\beta(s)$ | $=\widetilde{a}(s\sqrt{K})+s\sqrt{K}$ |
| $\widetilde{a}(s)$ | $=\frac{(1+\mu)^2+\psi_1(s)+\Lambda_2\psi_2(s)}{\min\{g_2(\zeta):|\zeta|\leq\Lambda_1+s\}}s$, if $s\geq\varepsilon$<br>$=\widetilde{k}s$, if $s\in[0,\varepsilon)$ |
| $P(s)$ | $P(s):=(2\Lambda_3\psi_1(s)+2\Lambda_2\Lambda_3\psi_2(s)+(\Lambda_4+\Lambda_3)s\psi_2(s))^2+\psi_3^2(s)$<br>$+1+\psi_4(s(1+\tau L(s)))(\Lambda_3+sL(s))^2+\Lambda_5\psi_3(s(1+\tau L(s)))$ |
| $\zeta_d(t)$ | $=\left(\tan(q_d(t)), \frac{\dot{q}_d(t)}{\cos^2(q_d(t))}\right)$ |
| $v_d(t-\tau)$ | $=\frac{\ddot{q}_d(t)+\frac{dF}{dq}(q_d(t))+H(\dot{q}_d(t))}{G(q_d(t), \dot{q}_d(t))}$ |

**Table 2:** Table of all functions involved the hybrid feedback law (2.6), (2.7), (2.8), (2.9). In all functions $\zeta\in\Re^2$, $x\in\Re^2$ and $s,t\geq 0$. The functions $\psi_i(s)$ ($i=1,...,4$) are arbitrary continuous, non-decreasing functions that satisfy
$\max\{|\nabla g_i(\zeta)|:|\zeta|\leq\Lambda_1+s\}\leq\psi_i(s)$, $i=1,2$, $\max\{|\nabla\widetilde{W}(x)|:|x|\leq\Lambda_1+s\}\leq\psi_3(s)$,
$\max\{|\nabla^2\widetilde{W}(x)|:|x|\leq\Lambda_1+s\}\leq\psi_4(s)$, for all $s\geq 0$. The functions $\theta_i\in K_\infty$ ($i=1,2$) and the constant $R_2\geq 0$ are selected in such a way that the inequality
$\theta_1(|x|)\leq\widetilde{W}(x)\leq R_2+\theta_2(|x|)$ holds for all $x\in\Re^2$.

## Acknowledgments

The work of M. Malisoff has been supported by NSF grant 1102348.

# Appendix

**Proof of Claim 1:** First we show that for each partition $\{T_i\}_{i=0}^{\infty}$ of $\Re_+$ with $\sup_{i \geq 0}(T_{i+1} - T_i) \leq r$, for each $(t_0, x_0) \in \Re_+ \times \Re^2$ and $u_0 \in L^{\infty}([-\tau, 0); \Re)$, the solution of (4.27), (4.28), (4.29) and (4.30) with initial condition $x(t_0) = x_0$, $\breve{T}_\tau(t_0)u = u_0$ is unique and exists for all $t \geq t_0$.

The solution of (4.27), (4.28), (4.29) and (4.30) is determined by the following process:

Initial Step: Given $x(t_0) = x_0$ and $\breve{T}_\tau(t_0)u = u_0$ we determine the solution $x(t)$ of (4.27) for $t \in [t_0, t_0 + \tau]$. Notice that the solution is unique. Inequality (3.9) implies the following estimate:

$$|x(t)| \leq a_\tau(|x_0| + \|u_0\|_\tau), \text{ for all } t \in [t_0, t_0 + \tau] \tag{A.1}$$

$i$-th Step: Given $x(t)$ for $t \in [t_0, t_0 + T_i + \tau]$ and $u(t)$ for $t \in [-\tau, t_0 + T_i)$ we determine $x(t)$ for $t \in [t_0, t_0 + T_{i+1} + \tau]$ and $u(t)$ for $t \in [-\tau, t_0 + T_{i+1})$. The solution $\xi(t)$ of (4.28) for $t \in [t_0 + T_i, t_0 + T_{i+1})$ with initial condition $\xi(t_0 + T_i) = z_{N_i}$ is unique (and is given by (2.8)). Inequality (4.14) implies:

$$V(\xi(t)) \leq V(\xi(t_0 + T_i)), \text{ for all } t \in [t_0 + T_i, t_0 + T_{i+1}) \tag{A.2}$$

We determine $u(t)$ for $t \in [t_0 + T_i, t_0 + T_{i+1})$ using the equation $u(t) = k(t, \xi(t))$. Notice that inequalities (4.15) and (4.17) in conjunction with (A.2) imply the following inequality for all $t \in [t_0 + T_i, t_0 + T_{i+1})$:

$$|u(t)| = |k(t, \xi(t))| \leq \tilde{a}(|\xi(t_0 + T_i)|\sqrt{K}) \tag{A.3}$$

Finally, we determine the solution $x(t)$ of (4.27) for $t \in [t_0, t_0 + T_{i+1} + \tau]$. Notice that the solution is unique. The fact that $T_{i+1} - T_i \leq r$ in conjunction with inequality (3.9) with $\tau$ replaced by $r > 0$ and inequality (A.3) implies this estimate:

$$|x(t)| \leq a_r(|x(t_0 + T_i + \tau)| + \tilde{a}(|\xi(t_0 + T_i)|\sqrt{K})), \text{ for all } t \in [t_0 + T_i + \tau, t_0 + T_{i+1} + \tau] \tag{A.4}$$

Next we evaluate the difference $\xi(t) - x(t + \tau)$ for $t \in [t_0 + T_i, t_0 + T_{i+1})$.

Exploiting (3.2) we get:

$$|\xi(t) - x(t+\tau)| =$$

$$= \left| \xi(t_0 + T_i) - x(t_0 + T_i + \tau) + \int_{t_0+T_i}^{t} (f(s+\tau, \xi(s), k(s+\tau, \xi(s))) - f(s+\tau, x(s+\tau), k(s+\tau, \xi(s)))) ds \right|$$

$$\leq |\xi(t_0 + T_i) - x(t_0 + T_i + \tau)| + \int_{t_0+T_i}^{t} L(|\xi(s)| + |x(s+\tau)| + |k(s+\tau, \xi(s))|) |\xi(s) - x(s+\tau)| ds$$

Using inequalities (4.15), (A.2), (A.3) and (A.4), in conjunction with the above inequality, we obtain:

$$|\xi(t) - x(t+\tau)| \leq |\xi(t_0 + T_i) - x(t_0 + T_i + \tau)|$$
$$+ L\left(|\xi(t_0 + T_i)|\sqrt{K} + \tilde{a}(|\xi(t_0 + T_i)|\sqrt{K}) + a_r(|x(t_0 + T_i + \tau)| + \tilde{a}(|\xi(t_0 + T_i)|\sqrt{K}))\right) \int_{t_0+T_i}^{t} |\xi(s) - x(s+\tau)| ds$$

Define $\varphi(s) := a_r(s) + s$. Using the Growall-Bellman lemma, the above inequality, formula $\beta(s) := \tilde{a}(s\sqrt{K}) + s\sqrt{K}$ from (4.20) and the fact that $T_{i+1} - T_i \leq r$, we get for all $t \in [t_0 + T_i, t_0 + T_{i+1})$:

$$|\xi(t) - x(t+\tau)| \leq |\xi(t_0 + T_i) - x(t_0 + T_i + \tau)| \exp\left(rL(\varphi(|x(t_0 + T_i + \tau)| + \beta(|\xi(t_0 + T_i)|)))\right) \quad (A.5)$$

Next we evaluate the quantity $\nabla V(x(t+\tau)) f(t+\tau, x(t+\tau), k(t+\tau, \xi(t)))$ for $t \in [t_0 + T_i, t_0 + T_{i+1})$. Using inequality (4.14) we get:
$$\nabla V(x(t+\tau)) f(t+\tau, x(t+\tau), k(t+\tau, \xi(t))) \leq -2\mu V(x(t+\tau)) +$$
$$\nabla V(x(t+\tau))(f(t+\tau, x(t+\tau), k(t+\tau, \xi(t))) - f(t+\tau, x(t+\tau), k(t+\tau, x(t+\tau))))$$

The following estimate follows from (4.17), (4.19) and the above inequality:

$$\nabla V(x(t+\tau)) f(t+\tau, x(t+\tau), k(t+\tau, \xi(t))) \leq -2\mu V(x(t+\tau)) +$$
$$2K|x(t+\tau)| M(|x(t+\tau)| + |\xi(t)|) |x(t+\tau) - \xi(t)|$$

Using the above inequality in conjunction with inequality (4.15), inequalities (A.2), (A.4) and definitions $\beta(s) := \tilde{a}(s\sqrt{K}) + s\sqrt{K}$ and $\varphi(s) := a_r(s) + s$, we get:

$$\nabla V(x(t+\tau)) f(t+\tau, x(t+\tau), k(t+\tau, \xi(t))) \leq -2\mu V(x(t+\tau)) +$$
$$2K(\varphi(|x(t_0 + T_i + \tau)| + \beta(|\xi(t_0 + T_i)|))) M(\varphi(|x(t_0 + T_i + \tau)| + \beta(|\xi(t_0 + T_i)|))) |x(t+\tau) - \xi(t)| \quad (A.6)$$

Combining inequalities (A.5), (A.6) and definition (4.20) we obtain the following for all $t \in [t_0 + T_i, t_0 + T_{i+1})$:

$$\nabla V(x(t+\tau)) f(t+\tau, x(t+\tau), k(t+\tau, \xi(t)))$$
$$\leq -2\mu V(x(t+\tau)) + D_r(|x(t_0 + T_i + \tau)| + \beta(|\xi(t_0 + T_i)|)) |x(t_0 + T_i + \tau) - \xi(t_0 + T_i)| \quad (A.7)$$

Since $\xi(t_0 + T_i) = z_{N_i}$ (recall (4.29)), it follows from (3.32) and (3.33) (applied with initial time $t_0 + T_i$), (2.2) and (2.4) that the following inequalities hold for all $i = 0, 1, 2, \ldots$:

$$|\xi(t_0 + T_i) - x(t_0 + T_i + \tau)| \leq R\left(|x(t_0 + T_i)| + \|\breve{T}_\tau(t_0 + T_i) u\|_\tau\right) \quad (A.8)$$

$$|\xi(t_0 + T_i)| \leq Q_\tau \left( |x(t_0 + T_i)| + \|\breve{T}_\tau(t_0 + T_i)u\|_\tau \right) \tag{A.9}$$

Since $|x(t_0 + T_i + \tau)| \leq a_\tau \left( |x(t_0 + T_i)| + \|\breve{T}_\tau(t_0 + T_i)u\|_\tau \right)$ (recall (3.9)), we obtain the following from (A.7), (A.8), (A.9) and definition (4.26) for all $t \in [t_0 + T_i, t_0 + T_{i+1})$:

$$\frac{d}{dt} V(x(t + \tau)) \leq -2\mu V(x(t + \tau)) + \gamma \tag{A.10}$$

Integrating the above differential inequality, we obtain for all $t \geq t_0$:

$$V(x(t + \tau)) \leq \exp(-2\mu(t - t_0)) V(x(t_0 + \tau)) + \frac{\gamma}{2\mu} \tag{A.11}$$

Combining (4.15), (A.1) and (A.11) we obtain inequality (4.34) with $\sigma(s,t) := 2K(a_\tau(s))^2 \exp(-2\mu(t - \tau))$ for all $t > \tau$ and $\sigma(s,t) := 2K(a_\tau(s))^2$ for all $t \in [0, \tau]$. The proof is complete. ◁

**Proof of Claim 2:** Let arbitrary partition $\{T_i\}_{i=0}^\infty$ of $\Re_+$ with $\sup_{i \geq 0}(T_{i+1} - T_i) \leq r$, $t_0 \geq 0$, $x_0 \in \Re^2$, $u_0 \in L^\infty([-\tau, 0); \Re^m)$ and consider the solution of (4.27), (4.28), (4.29) and (4.30) with (arbitrary) initial conditions $x(t_0) = x_0$ and $\breve{T}_\tau(t_0)u = u_0$. Inequalities (4.22) and (4.34) guarantee that there exists a unique smallest sampling time $t_0 + T_j$ such that $V(x(t_0 + T_j + \tau)) \leq \delta$, since $\frac{\gamma}{\mu} < \frac{\delta}{2}$.

Moreover, inequalities (A.10), (4.22) and (4.15) allow us to conclude that

$$|x(t)| \leq \sqrt{\delta} \text{ and } V(x(t)) \leq \delta, \text{ for all } t \geq t_0 + T_j + \tau \tag{A.12}$$

Using (A.8), definition (4.26), (4.22) and (A.12) we obtain that:

$$|\xi(t_0 + T_i)| \leq |\xi(t_0 + T_i) - x(t_0 + T_i + \tau)| + |x(t_0 + T_i + \tau)| \leq \gamma + \sqrt{\delta} \leq 2\sqrt{\delta}, \text{ for all } i \geq j \tag{A.13}$$

Using (A.2), (4.15) and (A.13), we get $|\xi(t)|^2 \leq V(\xi(t)) \leq 4\delta K$ for all $t \geq t_0 + T_j$, hence:

$$|\xi(t)| \leq 2\sqrt{K\delta} \tag{A.14}$$

for all $t \geq t_0 + T_j$.

Next we evaluate the difference $\xi(t) - x(t + \tau)$ for $t \geq t_0 + T_j$. Exploiting (3.2) for the vector field $f(t, x, u) := \begin{bmatrix} x_2 \\ \tilde{f}(t, x) + \tilde{g}(t, x)u \end{bmatrix}$ and inequalities (4.17), (4.18), (4.21), (A.12), (A.14) and definition (4.24), we get the following for all $i \geq j$ and $t \in [t_0 + T_i, t_0 + T_{i+1})$:

$$\left|\xi(t)-x(t+\tau)\right|$$
$$=\left|\xi(t_0+T_i)-x(t_0+T_i+\tau)+\int_{t_0+T_i}^{t}\bigl(f(s+\tau,\xi(s),k(s+\tau,\xi(s)))-f(s+\tau,x(s+\tau),k(s+\tau,\xi(s)))\bigr)ds\right|$$
$$\leq\left|\xi(t_0+T_i)-x(t_0+T_i+\tau)\right|+\widetilde{L}\int_{t_0+T_i}^{t}\left|\xi(s)-x(s+\tau)\right|ds$$

Using the Growall-Bellman lemma, the above inequality and the fact that $T_{i+1}-T_i \leq r$ imply that for all $i \geq j$ and $t \in [t_0+T_i, t_0+T_{i+1})$:

$$\left|\xi(t)-x(t+\tau)\right|\leq\left|\xi(t_0+T_i)-x(t_0+T_i+\tau)\right|\exp(r\widetilde{L}) \tag{A.15}$$

Next we evaluate the quantity $\nabla V(x(t+\tau))f(t+\tau,x(t+\tau),k(t+\tau,\xi(t)))$ for $t\in[t_0+T_i,t_0+T_{i+1})$. Using inequalities (4.14), (4.16), (4.18), (A.12), (4.21), (4.19), (A.14) and (A.15) and definition (4.23), we get the following for all $i\geq j$ and $t\in[t_0+T_i,t_0+T_{i+1})$:

$$\begin{aligned}
&\nabla V(x(t+\tau))f(t+\tau,x(t+\tau),k(t+\tau,\xi(t))) \\
&\leq -2\mu V(x(t+\tau))+2K|x(t+\tau)|\bigl\|f(t+\tau,x(t+\tau),k(t+\tau,\xi(t)))-f(t+\tau,x(t+\tau),k(t+\tau,x(t+\tau)))\bigr\| \\
&\leq -2\mu V(x(t+\tau))+2K|x(t+\tau)|\widetilde{g}(t,x(t+\tau))|k(t+\tau,\xi(t))-k(t+\tau,x(t+\tau))| \\
&\leq -2\mu V(x(t+\tau))+2K|x(t+\tau)|M\bigl(|x(t+\tau)|+|\xi(t)|\bigr)|\xi(t)-x(t+\tau)| \\
&\leq -2\mu V(x(t+\tau))+\phi|x(t+\tau)|\,|x(t_0+T_i+\tau)-\xi(t_0+T_i)|
\end{aligned} \tag{A.16}$$

Using (A.16) and the triangle inequality, we get the following for all $i\geq j$ and $t\in[t_0+T_i,t_0+T_{i+1})$:

$$\dot{V}(t+\tau)\leq -\mu V(t+\tau)+\frac{\phi^2}{4\mu}\left|x(t_0+T_i+\tau)-\xi(t_0+T_i)\right|^2 \tag{A.17}$$

where $V(t)=V(x(t))$. Using (A.8) and the fact that $R(s)\leq \widetilde{R}s$ for all $s\geq 0$, we get the following for all $i\geq j$ and $t\in[t_0+T_i,t_0+T_{i+1})$:

$$\dot{V}(t+\tau)\leq -\mu V(t+\tau)+\frac{\phi^2\widetilde{R}^2}{2\mu}\left|x(t_0+T_i)\right|^2+\frac{\phi^2\widetilde{R}^2}{2\mu}\left\|\widetilde{T}_\tau(t_0+T_i)u\right\|_\tau^2 \tag{A.18}$$

Let $\omega\in\left(0,\frac{\mu}{2}\right)$ be a positive constant sufficiently small such that

$$\widetilde{R}\widetilde{k}\sqrt{K}\exp(\omega(r+\tau))<1 \text{ and } \frac{\phi\widetilde{R}}{\sqrt{2\mu}}\frac{\exp(\omega(r+\tau))}{\sqrt{\mu-2\omega}}\left(1+\frac{\widetilde{k}\exp(\omega r)\sqrt{K}\bigl(\widetilde{R}+\exp(-\omega\tau)\bigr)}{1-\widetilde{R}\widetilde{k}\sqrt{K}\exp(\omega(r+\tau))}\right)<1 \tag{A.19}$$

The existence of $\omega\in\left(0,\frac{\mu}{2}\right)$ satisfying (A.19) is guaranteed by (4.25). Using (A.18) and the fact that $\sup_{i\geq 0}(T_{i+1}-T_i)\leq r$, we obtain the following for all $i\geq j$ and $t\in[t_0+T_i,t_0+T_{i+1})$:

$$\begin{aligned}
\dot{V}(t+\tau) &\leq -\mu V(t+\tau)+\frac{\phi^2\widetilde{R}^2}{2\mu}\exp(-2\omega t)\exp(2\omega r)\sup_{t_0+T_i\leq s\leq t}\bigl(\exp(2\omega s)|x(s)|^2\bigr) \\
&+\frac{\phi^2\widetilde{R}^2}{2\mu}\exp(-2\omega t)\exp(2\omega(r+\tau))\sup_{t_0+T_i-\tau\leq s\leq t}\bigl(\exp(2\omega s)|u(s)|^2\bigr)
\end{aligned} \tag{A.20}$$

The differential inequality (A.20) allows us to conclude that the following differential inequality holds for almost all $t \geq t_0 + T_j$:

$$\dot{V}(t+\tau) \leq -\mu V(t+\tau) + \frac{\phi^2 \widetilde{R}^2}{2\mu} \exp(-2\omega t)\exp(2\omega r) \sup_{t_0+T_j \leq s \leq t} \left(\exp(2\omega s)|x(s)|^2\right)$$
$$+ \frac{\phi^2 \widetilde{R}^2}{2\mu} \exp(-2\omega t)\exp(2\omega(r+\tau)) \sup_{t_0+T_j-\tau \leq s \leq t} \left(\exp(2\omega s)|u(s)|^2\right) \quad (A.21)$$

Multiplying (A.21) through by $\exp(\mu(t+\tau))$ and then integrating the result over $[t_0+T_j, t]$ for any $t \geq t_0 + T_j$ and using the fact that $2\omega < \mu$, we obtain the following for all $t \geq t_0 + T_j$:

$$V(t+\tau) \leq \exp(-2\omega(t-t_0-T_j))V(t_0+T_j+\tau) + \frac{\phi^2 \widetilde{R}^2}{2\mu} \frac{\exp(-2\omega t)}{\mu - 2\omega}\exp(2\omega r) \sup_{t_0+T_j \leq s \leq t} \left(\exp(2\omega s)|x(s)|^2\right)$$
$$+ \frac{\phi^2 \widetilde{R}^2}{2\mu} \frac{\exp(-2\omega t)}{\mu - 2\omega}\exp(2\omega(r+\tau)) \sup_{t_0+T_j-\tau \leq s \leq t} \left(\exp(2\omega s)|u(s)|^2\right) \quad (A.22)$$

Using our quadratic upper and lower bounds for $V$ from (4.15), we conclude from (A.22) that for all $t \geq t_0 + T_j$, we get:

$$|x(t+\tau)|\exp(\omega(t+\tau)) \leq \sqrt{K}\exp(\omega(t_0+T_j+\tau))|x(t_0+T_j+\tau)|$$
$$+ \frac{\phi \widetilde{R}}{\sqrt{2\mu}} \frac{\exp(\omega(r+\tau))}{\sqrt{\mu-2\omega}} \sup_{t_0+T_j \leq s \leq t} \left(\exp(\omega s)|x(s)|\right) \quad (A.23)$$
$$+ \frac{\phi \widetilde{R}}{\sqrt{2\mu}} \frac{\exp(\omega(r+2\tau))}{\sqrt{\mu-2\omega}} \sup_{t_0+T_j-\tau \leq s \leq t} \left(\exp(\omega s)|u(s)|\right)$$

Recall from (A.14) and our condition (4.21) on $\delta$ that $|\xi(t)| \leq \varepsilon$ for all $t \geq t_0 + T_j$. Hence, using our bounds (4.17)-(4.18) on $k$, (4.26), (A.2), our quadratic bounds (4.15) for $V$, (A.8), the fact that $R(s) \leq \widetilde{R}s$ for all $s \geq 0$, (A.14) and the triangle inequality, we obtain the following for all $i \geq j$ and $t \in [t_0+T_i, t_0+T_{i+1})$:

$$|u(t)| = |k(t+\tau, \xi(t))| \leq \widetilde{k}|\xi(t)| \leq \widetilde{k}\sqrt{K}|\xi(t_0+T_i)|$$
$$\leq \widetilde{k}\sqrt{K}|\xi(t_0+T_i) - x(t_0+T_i+\tau)| + \widetilde{k}\sqrt{K}|x(t_0+T_i+\tau)| \quad (A.24)$$
$$\leq \widetilde{R}\widetilde{k}\sqrt{K}|x(t_0+T_i)| + \widetilde{R}\widetilde{k}\sqrt{K}\|\widetilde{T}_\tau(t_0+T_i)u\|_\tau + \widetilde{k}\sqrt{K}|x(t_0+T_i+\tau)|$$

Inequality (A.24) in conjunction with the fact that $\sup_{i \geq 0}(T_{i+1}-T_i) \leq r$ implies the following for all $i \geq j$:

$$|u(t)|\exp(\omega t) \leq \widetilde{R}\widetilde{k}\sqrt{K}\exp(\omega r)|x(t_0+T_i)|\exp(\omega(t_0+T_i))$$
$$+ \widetilde{R}\widetilde{k}\sqrt{K}\exp(\omega(r+\tau)) \sup_{t_0+T_i-\tau \leq s < t_0+T_i} \left(\exp(\omega s)|u(s)|\right)$$
$$+ \widetilde{k}\sqrt{K}\exp(\omega(r-\tau))|x(t_0+T_i+\tau)|\exp(\omega(t_0+T_i+\tau))$$

The above inequality gives the following for all $t \geq t_0 + T_j$:

$$|u(t)|\exp(\omega t) \leq \widetilde{k}\exp(\omega r)\sqrt{K}\left(\widetilde{R}+\exp(-\omega\tau)\right)\sup_{t_0+T_j-\tau \leq s \leq t}\left(\exp(\omega(s+\tau))|x(s+\tau)|\right)$$
$$+ \widetilde{R}\widetilde{k}\sqrt{K}\exp(\omega(r+\tau))\sup_{t_0+T_j-\tau \leq s \leq t}\left(\exp(\omega s)|u(s)|\right) \quad (A.25)$$

Distinguishing the cases $\sup_{t_0+T_j-\tau \leq s \leq t}(\exp(\omega s)|u(s)|) = \sup_{t_0+T_j \leq s \leq t}(\exp(\omega s)|u(s)|)$ and $\sup_{t_0+T_j-\tau \leq s \leq t}(\exp(\omega s)|u(s)|) = \sup_{t_0+T_j-\tau \leq s < t_0+T_j}(\exp(\omega s)|u(s)|)$ we obtain the following from (A.25) for all $t \geq t_0+T_j$:

$$\sup_{t_0+T_j \leq s \leq t}(\exp(\omega s)|u(s)|) \leq \frac{\tilde{k}\exp(\omega r)\sqrt{K}(\tilde{R}+\exp(-\omega\tau))}{1-\tilde{R}\tilde{k}\sqrt{K}\exp(\omega(r+\tau))} \sup_{t_0+T_j-\tau \leq s \leq t}(\exp(\omega(s+\tau))|x(s+\tau)|)$$
$$+ \tilde{R}\tilde{k}\sqrt{K}\exp(\omega(r+\tau)) \sup_{t_0+T_j-\tau \leq s < t_0+T_j}(\exp(\omega s)|u(s)|) \quad (A.26)$$

Again, by distinguishing the cases $\sup_{t_0+T_j-\tau \leq s \leq t}(\exp(\omega s)|u(s)|) = \sup_{t_0+T_j \leq s \leq t}(\exp(\omega s)|u(s)|)$ and $\sup_{t_0+T_j-\tau \leq s \leq t}(\exp(\omega s)|u(s)|) = \sup_{t_0+T_j-\tau \leq s < t_0+T_j}(\exp(\omega s)|u(s)|)$, using the fact that $\tilde{R}\tilde{k}\sqrt{K}\exp(\omega(r+\tau)) < 1$ and combining (A.23) and (A.26), we get the following for all $t \geq t_0+T_j$:

$$|x(t+\tau)|\exp(\omega(t+\tau)) \leq \sqrt{K}\exp(\omega(t_0+T_j+\tau))|x(t_0+T_j+\tau)|$$
$$+ \frac{\phi\tilde{R}}{\sqrt{2\mu}}\frac{\exp(\omega(r+\tau))}{\sqrt{\mu-2\omega}}\left(1+\frac{\tilde{k}\exp(\omega r)\sqrt{K}(\tilde{R}+\exp(-\omega\tau))}{1-\tilde{R}\tilde{k}\sqrt{K}\exp(\omega(r+\tau))}\right)\sup_{t_0+T_j-\tau \leq s \leq t}(\exp(\omega(s+\tau))|x(s+\tau)|)$$
$$+ \frac{\phi\tilde{R}}{\sqrt{2\mu}}\frac{\exp(\omega(r+2\tau))}{\sqrt{\mu-2\omega}} \sup_{t_0+T_j-\tau \leq s < t_0+T_j}(\exp(\omega s)|u(s)|)$$

Distinguishing between the cases $\sup_{t_0+T_j-\tau \leq s \leq t}(\exp(\omega(s+\tau))|x(s+\tau)|) = \sup_{t_0+T_j-\tau \leq s \leq t_0+T_j}(\exp(\omega(s+\tau))|x(s+\tau)|)$ and $\sup_{t_0+T_j-\tau \leq s \leq t}(\exp(\omega(s+\tau))|x(s+\tau)|) = \sup_{t_0+T_j \leq s \leq t}(\exp(\omega(s+\tau))|x(s+\tau)|)$ and using the above inequality, we obtain the following for all $t \geq t_0+T_j$:

$$|x(t+\tau)|\exp(\omega(t+\tau)) \leq \frac{\sqrt{K}\exp(\omega(t_0+T_j+\tau))}{1-\lambda}|x(t_0+T_j+\tau)|$$
$$+ \frac{\phi\tilde{R}}{\sqrt{2\mu}}\frac{\exp(\omega(r+\tau))}{\sqrt{\mu-2\omega}}\left(1+\frac{\tilde{k}\exp(\omega r)\sqrt{K}(\tilde{R}+\exp(-\omega\tau))}{1-\tilde{R}\tilde{k}\sqrt{K}\exp(\omega(r+\tau))}\right)\sup_{t_0+T_j-\tau \leq s \leq t_0+T_j}(\exp(\omega(s+\tau))|x(s+\tau)|) \quad (A.27)$$
$$+ \frac{\phi\tilde{R}}{\sqrt{2\mu}}\frac{\exp(\omega(r+2\tau))}{(1-\lambda)\sqrt{\mu-2\omega}} \sup_{t_0+T_j-\tau \leq s < t_0+T_j}(\exp(\omega s)|u(s)|)$$

where $\lambda := \frac{\phi\tilde{R}}{\sqrt{2\mu}}\frac{\exp(\omega(r+\tau))}{\sqrt{\mu-2\omega}}\left(1+\frac{\tilde{k}\exp(\omega r)\sqrt{K}(\tilde{R}+\exp(-\omega\tau))}{1-\tilde{R}\tilde{k}\sqrt{K}\exp(\omega(r+\tau))}\right)$. Inequalities (A.26) and (A.27) imply that there exist positive constants $S_1$ and $S_2$ such that (4.35) and (4.36) hold.
The proof is complete. ◁

**Proof of Claim 3:** Let arbitrary partition $\{T_i\}_{i=0}^{\infty}$ of $\Re_+$ with $\sup_{i \geq 0}(T_{i+1}-T_i) \leq r$, $(t_0, x_0) \in \Re_+ \times \Re^2$ and $u_0 \in L^{\infty}([-\tau,0);\Re)$ be given and consider the solution of (4.27), (4.28), (4.29) and (4.30) with (arbitrary) initial conditions $x(t_0) = x_0$ and $\breve{T}_\tau(t_0)u = u_0$.

Define:

$$b(s) := \tilde{a}\left(s\sqrt{K}\right), \text{ for all } s \geq 0 \tag{A.28}$$

where $\tilde{a} \in K_\infty$ is from (4.17). Then $b \in K_\infty$. Moreover, notice that definitions (A.28) and (4.26) imply that

$$R(s) \leq \frac{1}{2} b^{-1}\left(\frac{s}{2}\right), \text{ for all } s \geq 0 \tag{A.29}$$

Furthermore, definition (A.28) and inequality (A.3) imply the following inequality for all $i \in Z_+$ and $t \in [t_0 + T_i, t_0 + T_{i+1})$:

$$|u(t)| \leq b(|\xi(t_0 + T_i)|) \tag{A.30}$$

Inequalities (4.15) and (4.34) imply the existence of a non-decreasing function $g: \Re_+ \to \Re_+$ such that:

$$|x(t)| \leq g(|x_0| + \|u_0\|_\tau), \text{ for all } t \geq t_0 \tag{A.31}$$

Combining (A.8), (A.29) and (A.31) we get the following for all $i \in Z_+$:

$$|\xi(t_0 + T_i) - x(t_0 + T_i + \tau)| \leq R\left(|x(t_0 + T_i)| + \sup_{t_0 + T_i - \tau \leq s < t_0 + T_i}(|u(s)|)\right)$$

$$\leq \frac{1}{2} b^{-1}\left(\frac{1}{2}|x(t_0 + T_i)| + \frac{1}{2}\sup_{t_0 + T_i - \tau \leq s < t_0 + T_i}(|u(s)|)\right)$$

$$\leq \max\left\{\frac{1}{2} b^{-1}(|x(t_0 + T_i)|), \frac{1}{2} b^{-1}\left(\sup_{t_0 + T_i - \tau \leq s < t_0 + T_i}(|u(s)|)\right)\right\}$$

$$\leq \max\left\{\frac{1}{2} b^{-1}(g(|x_0| + \|u_0\|_\tau)), \frac{1}{2} b^{-1}\left(\sup_{t_0 + T_i - \tau \leq s < t_0 + T_i}(|u(s)|)\right)\right\}$$

The above inequality in conjunction with (A.31) gives the following for all $i \in Z_+$:

$$|\xi(t_0 + T_i)| \leq |x(t_0 + T_i + \tau)| + \max\left\{\frac{1}{2} b^{-1}(g(|x_0| + \|u_0\|_\tau)), \frac{1}{2} b^{-1}\left(\sup_{t_0 + T_i - \tau \leq s < t_0 + T_i}(|u(s)|)\right)\right\}$$

$$\leq g(|x_0| + \|u_0\|_\tau) + \frac{1}{2}\max\left\{b^{-1}(g(|x_0| + \|u_0\|_\tau)), b^{-1}\left(\sup_{t_0 + T_i - \tau \leq s < t_0 + T_i}(|u(s)|)\right)\right\}$$

$$\leq \max\left\{2g(|x_0| + \|u_0\|_\tau), b^{-1}(g(|x_0| + \|u_0\|_\tau)), b^{-1}\left(\sup_{t_0 + T_i - \tau \leq s < t_0 + T_i}(|u(s)|)\right)\right\}$$

$$\leq \max\left\{2g(|x_0| + \|u_0\|_\tau), b^{-1}(g(|x_0| + \|u_0\|_\tau)), b^{-1}\left(\sup_{t_0 - \tau \leq s < t_0 + T_i}(|u(s)|)\right)\right\}$$

where we have used the inequalities $a_1 + a_2 \leq \max\{2a_1, 2a_2\}$ and $\max\{\lambda a_1, \lambda a_2\} = \lambda \max\{a_1, a_2\}$, which hold for all $a_i \in \Re_+$ ($i = 1, 2$) and $\lambda \geq 0$. Furthermore, using (A.30) and the above inequality, we obtain the following for all $i \in Z_+$:

$$\sup_{t_0 + T_i \leq s < t_0 + T_{i+1}} |u(s)| \leq \max\left\{\hat{g}(|x_0| + \|u_0\|_\tau), \sup_{t_0 - \tau \leq s < t_0 + T_i}(|u(s)|)\right\} \tag{A.32}$$

where $\hat{g}(s) := \max\{g(s), b(2g(s))\}$ for all $s \geq 0$, is a non-decreasing function. Define the sequence:

$$F_i := \sup_{t_0 - \tau \leq s < t_0 + T_i} (|u(s)|) \quad (A.33)$$

Notice that definition (A.33) and the fact that $\sup_{t_0 - \tau \leq s < t_0 + T_{i+1}} (|u(s)|) = \max\left\{ \sup_{t_0 + T_i \leq s < t_0 + T_{i+1}} (|u(s)|), \sup_{t_0 - \tau \leq s < t_0 + T_i} (|u(s)|) \right\}$ in conjunction with (A.32) imply the following inequality for all $i \in Z_+$:

$$F_{i+1} \leq \max\{\hat{g}(|x_0| + \|u_0\|_\tau), F_i\} \quad (A.34)$$

Inequality (A.34) in conjunction with the fact that $F_0 := \|u_0\|_\tau$ allow us to prove by induction that the following inequality holds for all $i \in Z_+$:

$$F_i \leq \max\{\hat{g}(|x_0| + \|u_0\|_\tau), \|u_0\|_\tau\} \quad (A.35)$$

Inequality (A.31) in conjunction with inequality (A.35) and definition (A.33) imply that estimate (4.37) holds with $S(s) := g(s) + \max\{\hat{g}(s), s\}$ for all $s \geq 0$. The proof is complete. ◁

**Construction of** $M : \Re_+ \to [1, +\infty)$ **satisfying (4.19):** Using (4.13) we obtain the following inequality for all $(t, x, \xi) \in \Re_+ \times \Re^2 \times \Re^2$

$$\tilde{g}(t,x)|k(t,x) - k(t,\xi)| \leq \left| (1+\mu^2)x_1 + 2\mu x_2 + \tilde{f}(t,x) - \frac{\tilde{g}(t,x)}{\tilde{g}(t,\xi)}\left((1+\mu^2)\xi_1 + 2\mu\xi_2 + \tilde{f}(t,\xi)\right) \right|$$

which after some rearranging gives:

$$\begin{aligned}&\tilde{g}(t,x)|k(t,x) - k(t,\xi)| \\ &\leq \left| (1+\mu^2)(x_1 - \xi_1) + 2\mu(x_2 - \xi_2) + \tilde{f}(t,x) - \tilde{f}(t,\xi) + \frac{\tilde{g}(t,\xi) - \tilde{g}(t,x)}{\tilde{g}(t,\xi)}\left((1+\mu^2)\xi_1 + 2\mu\xi_2 + \tilde{f}(t,\xi)\right) \right|\end{aligned} \quad (A.36)$$

Using (A.36), the facts that $|x_i - \xi_i| \leq |x - \xi|$ ($i = 1,2$) and the triangle inequality we obtain the following inequality for all $(t, x, \xi) \in \Re_+ \times \Re^2 \times \Re^2$:

$$\begin{aligned}&\tilde{g}(t,x)|k(t,x) - k(t,\xi)| \\ &\leq (1+\mu)^2|\xi - x| + \left|\tilde{f}(t,x) - \tilde{f}(t,\xi)\right| + \left|(1+\mu^2)\xi_1 + 2\mu\xi_2 + \tilde{f}(t,\xi)\right|\frac{|\tilde{g}(t,\xi) - \tilde{g}(t,x)|}{\tilde{g}(t,\xi)}\end{aligned} \quad (A.37)$$

Using (3.2), (3.3) and (4.6) we get $\left|\tilde{f}(t,x) - \tilde{f}(t,\xi)\right| \leq L(|x| + |\xi|)|\xi - x|$ and $\left|\tilde{f}(t,\xi)\right| \leq L(|\xi|)|\xi|$ for all $(t, x, \xi) \in \Re_+ \times \Re^2 \times \Re^2$. Combining the previous inequalities with (A.37) and using the facts that $|\xi_i| \leq |\xi|$ ($i = 1,2$) and the triangle inequality, we get the following inequality for all $(t, x, \xi) \in \Re_+ \times \Re^2 \times \Re^2$:

$$\begin{aligned}&\tilde{g}(t,x)|k(t,x) - k(t,\xi)| \\ &\leq (1+\mu)^2|\xi - x| + L(|x| + |\xi|)|\xi - x| + \left((1+\mu)^2 + L(|\xi|)\right)|\xi|\frac{|\tilde{g}(t,\xi) - \tilde{g}(t,x)|}{\tilde{g}(t,\xi)}\end{aligned} \quad (A.38)$$

Let $\psi_2(s)$ $(i=1,...,4)$ be a continuous, non-decreasing function that satisfies $\max\{|\nabla g_2(\zeta)|:|\zeta|\leq \Lambda_1+s\}\leq \psi_2(s)$, for all $s\geq 0$, where $\Lambda_1 := \sup_{t\geq 0}|\zeta_d(t)|$. Using definition (4.6) in conjunction with (A.37), we get the following inequality for all $(t,x,\xi)\in \Re_+\times\Re^2\times\Re^2$:

$$\tilde{g}(t,x)|k(t,x)-k(t,\xi)|$$
$$\leq (1+\mu)^2|\xi-x|+L(|x|+|\xi|)|\xi-x|+\left((1+\mu)^2+L(|\xi|)\right)|\xi|\frac{\psi_2(|x|+|\xi|)}{\min\{g_2(\zeta):|\zeta|\leq \Lambda_1+|\xi|\}}|x-\xi| \quad (A.39)$$

Finally, using the facts that $L(|\xi|)\leq L(|\xi|+|x|)$, $|\xi|\leq|\xi|+|x|$, $\min\{g_2(\zeta):|\zeta|\leq \Lambda_1+|\xi|+|x|\}\leq \min\{g_2(\zeta):|\zeta|\leq \Lambda_1+|\xi|\}$, in conjunction with (A.39), we are in a position to conclude that inequality (4.19) holds with

$$M(s):=\left((1+\mu)^2+L(s)\right)\left(1+\frac{s\psi_2(s)}{\min\{g_2(\zeta):|\zeta|\leq \Lambda_1+s\}}\right), \text{ for all } s\geq 0$$

The construction is complete. ◁

**Construction of the function** $P:\Re_+\to\Re_+$ **satisfying (3.5), (3.6) and (3.7):** Definition (4.6) implies that the following inequality holds for all $(s,t,x,u)\in\Re_+\times\Re_+\times\Re^2\times\Re$:

$$|f(s,x,u)-f(t,x,u)|\leq |s-t|\sup_{l\geq 0}\left|\frac{\partial}{\partial l}\tilde{f}(l,x)\right|+|s-t||u|\sup_{l\geq 0}\left|\frac{\partial}{\partial l}\tilde{g}(l,x)\right| \quad (A.40)$$

Let $\psi_i(s)$ $(i=1,2)$ be continuous, non-decreasing functions that satisfy $\max\{|\nabla g_i(\zeta)|:|\zeta|\leq \Lambda_1+s\}\leq \psi_i(s)$, for all $s\geq 0$, where $\Lambda_1 := \sup_{t\geq 0}|\zeta_d(t)|$. Inequality (A.40) in conjunction with the previous inequalities and definition (4.6), gives the following inequality for all $(s,t,x,u)\in\Re_+\times\Re_+\times\Re^2\times\Re$:

$$|f(s,x,u)-f(t,x,u)|\leq |s-t|\left(2\psi_1(|x|)\sup_{l\geq 0}|\dot{\zeta}_d(l)|+2\psi_2(|x|)\sup_{l\geq -\tau}|v_d(l)|\sup_{l\geq 0}|\dot{\zeta}_d(l)|+|x|\psi_2(|x|)\sup_{l\geq -\tau}|\dot{v}_d(l)|\right)$$
$$+|s-t||u|\psi_2(|x|)\sup_{l\geq 0}|\dot{\zeta}_d(l)| \quad (A.41)$$

Therefore, inequality (A.41) implies that inequality (3.7) holds provided that the following inequality holds for all $s\geq 0$:

$$2\Lambda_3\psi_1(s)+2\Lambda_2\Lambda_3\psi_2(s)+(\Lambda_4+\Lambda_3)s\psi_2(s)\leq \sqrt{P(s)} \quad (A.42)$$

where $\Lambda_2 := \sup_{t\geq -\tau}|v_d(t)|$, $\Lambda_3 := \sup_{t\geq 0}|\dot{\zeta}_d(t)|$ and $\Lambda_4 := \sup_{t\geq -\tau}|\dot{v}_d(t)|$.

Let $\tilde{W}(x)$ be the smooth function defined by $\tilde{W}(x):=1+\frac{1}{2}\left(\frac{x_2}{1+x_1^2}\right)^2+F\left(\tan^{-1}(x_1)\right)$ for all $x\in\Re^2$. Let $\psi_i(s)$ $(i=3,4)$ be ontinuous, non-decreasing functions that satisfy $\max\{|\nabla\tilde{W}(x)|:|x|\leq \Lambda_1+s\}\leq \psi_3(s)$, and

$\max\left\{\left|\nabla^2 \widetilde{W}(x)\right|:|x|\leq \Lambda_1+s\right\}\leq \psi_4(s)$, for all $s\geq 0$. Using the fact that $W(t,x)=\widetilde{W}(\zeta_d(t)+x)$ and the previous inequalities, we get the following inequality for all $(t,x)\in \mathfrak{R}_+ \times \mathfrak{R}^2$:

$$\left|\frac{\partial W}{\partial x}(t,x)\right|\leq \psi_3(|x|) \tag{A.43}$$

Therefore, inequality (A.43) implies that inequality (3.6) holds provided that the following inequality holds for all $s\geq 0$:

$$\psi_3(s)\leq \sqrt{P(s)} \tag{A.44}$$

Finally, we remark that the following inequality holds for all $s\geq 0$:

$$1+\sup\left\{\left|\frac{\partial^2 W}{\partial t^2}(t,\xi)\right|+2sL(s)\left|\frac{\partial^2 W}{\partial t\partial x}(t,\xi)\right|+s^2 L^2(s)\left|\frac{\partial^2 W}{\partial x^2}(t,\xi)\right|:|\xi|\leq s(1+\tau L(s)), t\geq 0\right\}$$
$$\leq 1+\sup\left\{\left|\frac{\partial^2 W}{\partial t^2}(t,\xi)\right|:|\xi|\leq s(1+\tau L(s)), t\geq 0\right\}+2sL(s)\sup\left\{\left|\frac{\partial^2 W}{\partial t\partial x}(t,\xi)\right|:|\xi|\leq s(1+\tau L(s)), t\geq 0\right\} \tag{A.45}$$
$$+s^2 L^2(s)\sup\left\{\left|\frac{\partial^2 W}{\partial x^2}(t,\xi)\right|:|\xi|\leq s(1+\tau L(s)), t\geq 0\right\}$$

Using the facts that $W(t,x)=\widetilde{W}(\zeta_d(t)+x)$, $\max\left\{\left|\nabla \widetilde{W}(x)\right|:|x|\leq \Lambda_1+s\right\}\leq \psi_3(s)$ and $\max\left\{\left|\nabla^2 \widetilde{W}(x)\right|:|x|\leq \Lambda_1+s\right\}\leq \psi_4(s)$, we conclude from (A.45) that the following inequality holds for all $s\geq 0$:

$$1+\sup\left\{\left|\frac{\partial^2 W}{\partial t^2}(t,\xi)\right|+2sL(s)\left|\frac{\partial^2 W}{\partial t\partial x}(t,\xi)\right|+s^2 L^2(s)\left|\frac{\partial^2 W}{\partial x^2}(t,\xi)\right|:|\xi|\leq s(1+\tau L(s)), t\geq 0\right\}$$
$$\leq 1+\psi_4(s(1+\tau L(s)))\left(\sup_{l\geq 0}|\dot\zeta_d(l)|\right)^2+\psi_3(s(1+\tau L(s)))\sup_{l\geq 0}|\ddot\zeta_d(l)|+2sL(s)\psi_4(s(1+\tau L(s)))\sup_{l\geq 0}|\dot\zeta_d(l)| \tag{A.46}$$
$$+s^2 L^2(s)\psi_4(s(1+\tau L(s)))$$

Therefore, inequality (A.46) implies that inequality (3.5) holds provided that the following inequality holds for all $s\geq 0$:

$$1+\psi_4(s(1+\tau L(s)))(\Lambda_3+sL(s))^2+\Lambda_5\psi_3(s(1+\tau L(s)))\leq P(s) \tag{A.47}$$

where $\Lambda_5:=\sup_{t\geq 0}|\ddot\zeta_d(t)|$ and $\Lambda_3:=\sup_{t\geq 0}|\dot\zeta_d(t)|$. Combining (A.42), (A.44) and (A.47) we are in a position to conclude that inequalities (3.5), (3.6) and (3.7) hold with

$$P(s):=(2\Lambda_3\psi_1(s)+2\Lambda_2\Lambda_3\psi_2(s)+(\Lambda_4+\Lambda_3)s\psi_2(s))^2+\psi_3^2(s)$$
$$+1+\psi_4(s(1+\tau L(s)))(\Lambda_3+sL(s))^2+\Lambda_5\psi_3(s(1+\tau L(s)))$$, for all $s\geq 0$

The construction is complete. ◁